\documentclass[12pt]{article}
\usepackage{mathrsfs}
\usepackage{epic,eepic,epsf,epsfig}
\usepackage{amsfonts,srcltx,mathrsfs}
\usepackage{amsmath}
\usepackage{amssymb}
\usepackage{amsbsy}
\usepackage{float}
\usepackage{graphicx}
\usepackage{amsfonts}
\usepackage{color}

\usepackage{url}
\newtheorem{lem}{Lemma}[section]
\newtheorem{theorem}[lem]{Theorem}

\newtheorem{cor}[lem]{Corollary}

\newtheorem{case}{Case}
\newtheorem{prop}[lem]{Proposition}

\newtheorem{conjecture}[lem]{Conjecture}

\def\a{\alpha}   \def\d{\delta} 
\def\r{\rho}

\def\olg{\overline G}

 \def\lg{\langle} \def\rg{\rangle}
\def\nd{\mathrel{\bigm|\kern-.7em/}}

\def\f{\noindent}

\def\M{\hbox{\mathcal M}}

\def\Aut{\hbox{\rm Aut}}

\def\Aut{\hbox{\rm Aut}}

\def\M{\hbox{\rm M}}

\def\demo{\f {\bf Proof.}\hskip10pt}

\newcommand{\qed}{\mbox{\raisebox{0.7ex}{\fbox{}}} \vspace{4truemm}}
\def\mz{{\mathbb Z}}
\def\M{\mathcal {M}}

\begin{document}

\title{Regular maps of order $2$-powers}

\author{Dong-Dong Hou$^a$, Yan-Quan Feng$^a$, Young Soo Kwon$^b$\\
{\small $^a$Department of Mathematics, Beijing Jiaotong University, Beijing,
100044, P.R. China}\\
{\small $^b$Mathematics, Yeungnam University, Kyongsan 712-749, Republic of Korea}}

\date{}
\maketitle

\footnotetext{
E-mails: holderhandsome$@$bjtu.edu.cn, yqfeng$@$bjtu.edu.cn, ysookwon$@$ynu.ac.kr}

\begin{abstract}

In this paper, we consider the possible types of regular maps of order $2^n$, where the order of a regular map is the order of  automorphism group of the map.
 For $n \le 11$, M. Conder classified all regular maps of order $2^n$. It is easy to classify  regular maps of order $2^n$ whose valency or covalency is $2$ or $2^{n-1}$.
 So we assume that $n \geq 12$ and $2\leq s,t\leq n-2$ with $s\leq t$ to consider regular maps of order $2^n$ with type $\{2^s, 2^t\}$. We show that for $s+t\leq n$ or for $s+t>n$ with $s=t$, there exists a regular map of order $2^n$ with type $\{2^s, 2^t\}$, and furthermore, we classify regular maps of order $2^n$ with types $\{2^{n-2},2^{n-2}\}$ and $\{2^{n-3},2^{n-3}\}$.
 We conjecture that, if $s+t>n$ with $s<t$, then there is no regular map of order $2^n$ with type $\{2^s, 2^t\}$, and we confirm the conjecture for $t=n-2$ and $n-3$.
 
\bigskip
\f {\bf Keywords:} Regular map, $2$-group, automorphism group.
\end{abstract}

\section{Introduction}
A \emph{map} $\M$ is a 2-cell embedding of a connected graph into a
closed surface.  The embedded graph is called the \emph{underlying
graph} of the map. A map is called orientable or nonorientable
according to whether the supporting surface is orientable or
non-orientable.  A map automorphism is a permutation of flags
(mutually  incident vertex-edge-face triples) that preserves their
relations of having a vertex, edge or face in common, namely, it
induces an automorphism of its underlying graph which extends to a
self-homeomorphism of its supporting surface. The automorphism group
of a map always acts freely on flags. If this action is transitive
as well, the map is said to be \emph{regular}. Obviously, if $\M$ is
regular, then all vertices of $\M$ have the same valency, say, $l$,
and all face boundary walks have the same length, say $m$. In this
case, $\M$ is said to have type $\{m, l\}$. If the order of automorphism group of a regular map is $n$,
we say that the regular map has order $n$, namely, the order of regular map equals to the number of flags in the map.
 The study of regular
maps has a long and rich history and was progressed substantially by
Brahana~\cite{H1927}, and Bryant and Singerman~\cite{RD1985}. For
more development of the theory of maps, we refer the readers to
~\cite{DCJ2014, CMa2015, NG2013}.

There are three different approaches to classify regular maps on surfaces: Classification of regular maps by underlying graphs~\cite{ARJM1999},
by map automorphism groups~\cite{G1994, C1969}, and by supporting surfaces~\cite{ARJ2005,MRJ2012,MP2001}.
In this paper, we concentrate on classifications of regular maps by automorphism groups.

Malni\v{c},  Nedela and \v{S}koviera~\cite{ARM2012} proved that each
regular map with a nilpotent automorphism group (which is called a
nilpotent regular map) can be uniquely decomposed into a direct
product of two regular maps: the automorphism group of one is a
$2$-group and the other map is a semistar of odd valency. This
implies that the automorphism group of a nilpotent regular map with
a simple underlying graph is a 2-group.
 So the classification of regular maps whose automorphism groups are
 2-groups is important to understand nilpotent regular maps. A regular map is
called a \emph{regular $2$-map} if its automorphism group is a
$2$-group. In~\cite{ARM2012}, Malni\v{c} \emph{et al} gave a
complete classification of nilpotent regular maps of nilpotency class
$2$. It is proved in~\cite{CDNS} that given the class, there are
finitely many simple regular $2$-maps. However, for the regular
$2$-maps with multiple edges and given class, it is possible to list
it by a computer and the classification of those seems to be
difficult. In~\cite{BDR}, regular 2-maps of class $3$ has been
classified. Furthermore, Hu \emph{et al.}~\cite{HW} classified
regular 2-maps for maximal class by using the classification of
2-groups with a cyclic maximal subgroup.

There are many papers~\cite{CHNS2012, GW, JNS, Si, Vi1983} showing  that for any given positive integers $m$ and $l$ satisfying $1/m + 1/l \le 1/2$, there exist infinitely many regular maps of type $\{m, l\}$ by constructive or non-constructive way. But, those papers deal with only existence of regular maps of given type without considering the order of automorphism group. In this paper, we consider the possible types of regular maps whose automorphism groups are 2-groups of order $2^n$ for a given integer $n$.
For $n \leq 11$, all such regular maps are listed in
\cite{atles}. It is easy to classify  regular maps of order $2^n$ whose valency or covalency is $2$ or $2^{n-1}$.
 So we assume that $n \geq 12$ and $2\leq s,t\leq n-2$ with $s\leq t$ to consider regular maps of order $2^n$ with type $\{2^s, 2^t\}$. We prove that for $s+t\leq n$, or $s+t>n$ and $s=t$, there exists a regular map of order $2^n$ with type $\{2^s, 2^t\}$, and for $s+t>n$, there is no regular map of order $2^n$ with type $\{2^s, 2^t\}$ when $t=n-2$ or $n-3$.  Furthermore, we classify regular maps of order $2^n$ with types $\{2^{n-2},2^{n-2}\}$ and $\{2^{n-3},2^{n-3}\}$.

This paper is organized as follows. In Section~\ref{Preliminaries},
we give some background and properties of regular maps. In Section~3 we consider the existence of a regular map $\M$ of order $2^n$ with given types, and in Section~4,  classifications of all regular maps of order $2^n$ with types $\{2^{n-2}, 2^{n-2}\}$ and $\{2^{n-3}, 2^{n-3}\}$ are given.

\section{Preliminaries}\label{Preliminaries}
In this section we present basic facts about regular maps. We
introduce regular maps by  starting from their automorphism groups
that are known to be quotients of extended triangle
groups~\cite{RD1985}. In all the forthcoming group presentations we
will assume that the listed exponents are the true orders of the
corresponding elements.

A finite regular map $\M$ can in this way be identified with a
(partial) three-generator presentation of a finite group $G$,
isomorphic to the automorphism group $\Aut(\M)$ of $\M$, of the form
$$\quad~~~ G=\lg \r_0, \r_1 ,\r_2\ |\ \r_0^2 = \r_1^2 = \r_2^2 = (\r_0\r_1)^m = (\r_1\r_2)^l = (\r_0\r_2)^2 = \cdots =1 \rg, \quad~~~~~~~~~~~~~~ (  *)$$
where dots indicate possible presence of additional relations.
In~\cite{LCH}, Li and \v{S}ir\'{a}\v{n} considered regular maps
whose automorphism groups correspond to the case that at least one
of $\r_0, \r_1, \r_2$ is the identity. So we just consider the case
that none of $\r_0, \r_1, \r_2$ is equal to $1$. We will construct a
regular map $\M$ such that $\Aut(\M)\cong G$.

Take topological triangles as many as the order of $G$. These
triangles are the \emph{flags} of our map to be constructed. Each
flag is labeled by an element of $G$, whereby distinct flags have
distinct labels. By this way, the set of flags is identified with
the set of elements of $G$. Moreover, the three sides of each flag
are displayed in blue, black and red lines. In each flag, the
blue, black and red sides are labeled by $\r_0$, $\r_1$ and
$\r_2$, respectively.

\vskip -1cm
\begin{figure}[H]
  \centering
  \includegraphics[width=14cm]{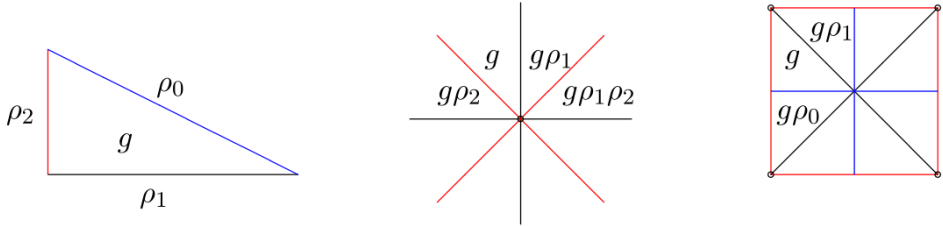}\\ \vskip -0.5cm
    \caption{Constructing a map from flags corresponding a group}\label{1}
\end{figure}

Next, for each $g \in G$ and each $w \in \{\r_0, \r_1, \r_2\}$, we
take the flags $g$ and $gw$ and identify their two sides whose
labels are  $w$. Applying this identification procedure with all
flags, one can form a closed surface $\mathcal{S}$. The union of all
red segments determines the underlying graph and its 2-cell
embedding on $\mathcal{S}$ constitutes our map $\M$. By $\M=\M(G;
\r_0, \r_1, \r_2)$ we denote the map constructed in this way and
call $\M$ the map associated with the presentation $(*)$ of $G$.
Considering this construction of map related to the group $G$, one
may identify vertices, edges and faces of the map $\M$ with the left
cosets of the subgroups $\lg \r_1, \r_2\rg$, $\lg \r_0, \r_2\rg$ and
$\lg \r_0, \r_1\rg$, respectively, and  their mutual incidence is
determined by non-empty intersection.

The group $G$ has two natural actions on the flag set of $\M$ (that
is, on itself), namely, by left and right multiplication.  The right
multiplication by the generators $\r_0, \r_1$, and $\r_2$ applied to
any flag is called the \emph{longitudinal reflection}, the
\emph{corner reflection}, and the \emph{transversal reflection} of
the flag. Regarding the action of $G$ by left multiplication, note
that if two flags $h_1, h_2 \in G$ are related by one of the above
reflections, then for any $g \in G$ the flags $gh_1$ and $gh_2$ are
related by  the same reflection. In this sense the left
multiplication corresponds to an automorphism of $M$. This enables
us to identify the (full) automorphism group $Aut(\M)$ of the map
$\M=\M(G; \r_0, \r_1, \r_2)$ with the group $G$ and its left action
on itself. In particular, regular maps of \emph{type} $\{m, l\}$
(that is, of face length $m$ and vertex valence $l$) can be
identified with presentations of finite $3$-generator groups as in
$(*)$. This allows one to translate the entire theory of regular
maps into a purely group-theoretical language. It is
known~\cite{RD1985, LCH} that the supporting surface of $\M$ is
orientable if and only if the rotation subgroup $\lg \r_0\r_1, \r_1\r_2\rg$ is a subgroup
of index $2$ in $G$. The following proposition shows when two given
regular maps $\M^{(1)}={\rm Map}(G^{(1)}, \{\r_0^{(1)}, \r_1^{(1)},
\r_2^{(1)}\})$ and $\M^{(2)}={\rm Map}(G^{(2)}, \{\r_0^{(2)},
\r_1^{(2)}, \r_2^{(2)}\})$ are isomorphic.


\begin{prop}\label{isomorphic}
$\M^{(1)}$ and $\M^{(2)}$ are isomorphic if and only if there is a group isomorphism $f: G^{(1)} \mapsto G^{(2)}$ such that $f(\r_0^{(1)})=\r_0^{(2)}$, $f(\r_1^{(1)})=\r_1^{(2)}$ and $f(\r_2^{(1)})=\r_2^{(2)}$.
\end{prop}

For a group $G$ with $x, y \in G$, denote by $[x, y]$ the \emph{commutators} $x^{-1}y^{-1}xy$ of $x$ and $y$, and the \emph{derived group} $G'$ of $G$ is the subgroup generated by all commutators $[x, y]$ for any $x, y\in G$.
The following  proposition is a basic property of commutators.

\begin{prop}\label{commutator}
Let $G$ be a group. Then, for any $x, y, z \in G$, $[xy, z]=[x, z]^y[y, z]$ and $[x, yz]=[x, z][x, y]^z$.
\end{prop}

With Proposition~\ref{commutator}, it is easy to prove that if $G$ is generated by a subset $M$, then $G'$ is generated by all conjugates in $G$ of elements $[x_i, x_j]$ with $x_i, x_j \in M$ (see ~\cite[Hilfsatz \uppercase\expandafter{\romannumeral3}.1.11]{GroupBooks}). Then it is easy to obtain the following result.

\begin{prop}{\rm ~\cite[Lemma 4.2]{HFL}}\label{derived}
Let $G=\lg \r_0, \r_1, \r_2 \ |\ \r_0^2, \r_1^2, \r_2^2, (\r_0\r_2)^2\rg$. Then $G'=\lg [\r_0, \r_1], [\r_1, \r_2], [\r_0, \r_1]^{\r_2}\rg$.
\end{prop}

The {\em Frattini subgroup}, denoted by $\Phi(G)$, of a finite group $G$ is defined to be the intersection of all maximal subgroups of $G$, and $\Phi(G)$ can be characterised as the set of non-generators of $G$.

\begin{prop}{\rm ~\cite[Theorem 9.2(a)]{2GroupBookss}}\label{frattini}
Let $G$ be a finite group and $S \subset G$. Then $G=\lg S \rg$ if and only if $G= \lg S, \Phi(G)\rg$.
\end{prop}

Let $G$ be a finite $p$-group for a prime $p$,  and set $\mho_1(G) = \lg g^p\ |\ g \in G\rg$. The following theorem is the well-known Burnside Basis Theorem.

\begin{prop}{\rm ~\cite[Theorem 1.12]{GroupBookss}}\label{burnside}
Let $G$ be a $p$-group and $|G: \Phi(G)| = p^d$.
\begin{itemize}
\item [(1)] $G/\Phi(G) \cong \mz_p^d$. Moreover, if $N \lhd G$ and $G/N$ is elementary abelian, then $\Phi(G) \leq N$.

\item [(2)] Every minimal generating set of $G$ contains exactly $d$ elements.

\item [(3)] $\Phi(G) = G' \mho_1(G)$. In particular, if $p=2$, then $\Phi(G) = \mho_1(G)$.
\end{itemize}
\end{prop}

The following is a classification of $2$-groups having a maximal cyclic subgroup.

\begin{prop}{\rm ~\cite[Theorem 1.2]{GroupBookss}}\label{maximal}
Let $G$ be a group of order $2^n$. If $G$ has a cyclic
subgroup of order $2^{n-1}$, then $G$ lies in one of the following six classes:

\begin{itemize}
  \item [{\rm (1)}] $G=\langle a\  | \ \ a^{2^{n}}=1\rangle, \ n\geq
  1$, the cyclic group.
\item [{\rm (2)}]
$G=\langle a,b\ |\ a^{2^{n-1}}=1,\  b^{2}=1,\  a^{b}=a\rangle,\
n\geq 2$, the abelian group.
\item [{\rm (3)}]$G=\langle a, b\ |\ a^{2^{n-1}}=1,\  b^{2}=1,\
 a^{b}=a^{-1}\rangle,\  n\geq 3$, the dihedral group. All elements in $G \backslash \lg a\rg$ are involutions.

\item [{\rm (4)}]
$G=\langle a,b\ |\ a^{2^{n-1}}=1,\  b^{2}=a^{2^{n-2}},\
a^{b}=a^{-1}\rangle,\  n\geq 3$, the generalized
quaternion group. The group G contains exactly one involution, and all elements in
$G \backslash \lg a \rg$ have the same order $4$.
\item [{\rm (5)}]
$G=\langle a, b\ |\ a^{2^{n-1}}=1, b^{2}=1,\
a^{b}=a^{1+2^{n-2}}\rangle,\  n\geq 4.$
\item [{\rm (6)}] $G=\langle
a, b\ |\ a^{2^{n-1}}=1, b^{2}=1,\  a^{b}=a^{-1+2^{n-2}}\rangle,\
n\geq 4$, the
semidihedral group.
\end{itemize}
\end{prop}

For a group $G$ and $g\in G$, we denote by $o(g)$ the order of $g$ in $G$. From Proposition~\ref{maximal}, we have the following result.

\begin{cor}\label{maximalc}
Let $n\geq 4$ and let $G$ be a group of order $2^n$ with a cyclic subgroup of order $2^{n-1}$. Then $u^v=u$ or $u^{1+2^{n-2}}$ for any $u, v \in G$ with $o(u)=o(v)=2^{n-1}$.
\end{cor}

\demo Let $G(i)$ be the groups defined in Proposition~\ref{maximal}~(i) for each $i \in
\{1, 2, 3, 4, 5, 6\}$.  Note that $G(1)$ and $G(2)$ are abelian. Since $n \geq 4$, it is easy to see that $G(3)$ and $G(4)$ has a unique cyclic subgroup of order $2^{n-1}$. For $G(6)$, we have $(a^jb)^2=a^j(a^j)^b=a^j(a^b)^j=a^{j\cdot 2^{n-2}}$ for $j \in \mz_{2^{n-1}}$, that is, $o(a^jb)=2$ or $4$, and hence $G(6)$ also has a unique cyclic subgroup of order $2^{n-1}$. Thus $\lg u\rg=\lg v\rg$ for $G(3)$, $G(4)$ and $G(6)$. It follows $u^v=u$ for $G=G(i)$ with $i=1,2,3,4$ or $6$.

Let $G=G(5)$. It is easy to see $G'=\lg a^{2^{n-2}}\rg\cong \mz_2$. Assume $u^v\not=u$. Then $\lg u\rg\not=\lg v\rg$ and $[u,v]=a^{2^{n-2}}$. Since $o(u)=o(v)=2^{n-1}$, we have
$G=\lg u,v\rg$,  which implies $\lg u\rg\cap \lg v\rg=\lg u^2\rg=\lg v^2\rg$. In particular, $\lg u^2\rg=\lg v^2\rg=\lg a^2\rg$, forcing $a^{2^{n-2}}=u^{2^{n-2}}$. It follows $u^v=u[u,v]=ua^{2^{n-2}}=u^{1+2^{n-2}}$, as required.
\hfill\qed

The following proposition gives the automorphism group of some regular maps of order $2^n$ for each $n \geq 10$.

\begin{prop}{\rm \cite[Theorem 1.2 and Theorem 1.3(1)]{HFL}}\label{$2^s,2^t$}
Let $n \geq 10$, $s, t \geq 2$ and $n-s-t \geq 1$. Set $R(\r_0, \r_1, \r_2)= \{\r_0^2, \r_1^2, \r_2^2, (\r_0\r_1)^{2^{s}}, (\r_1\r_2)^{2^{t}}, (\r_0\r_2)^2$, $[(\r_0\r_1)^4, \r_2], [\r_0,(\r_1\r_2)^4]\}$ and define
$$H=\left\{
\begin{array}{ll}
\lg \r_0, \r_1, \r_2 \ |\ R(\r_0, \r_1, \r_2), [(\r_0\r_1)^2, \r_2]^{2^{\frac{n-s-t-1}{2}}}\rg, & n-s-t\mbox{ odd }\\
\lg \r_0, \r_1, \r_2 \ |\ R(\r_0, \r_1, \r_2), [(\r_0\r_1)^2, (\r_1\r_2)^2]^{2^{\frac{n-s-t-2}{2}}}\rg, & n-s-t\mbox{ even. }
\end{array}
\right.$$ and $L=\lg \r_0, \r_1, \r_2 \ |\ \r_0^2, \r_1^2, \r_2^2, (\r_0\r_1)^{2^{2}}, (\r_1\r_2)^{2^{n-3}}, (\r_0\r_2)^2, [(\r_0\r_1)^2,\r_2](\r_1\r_2)^{2^{n-4}} \rg$.
Then $|H|=|L|=2^{n}$ and
the listed exponents are the true orders of the corresponding elements.

\end{prop}

\section{Existence of regular maps of order $2^n$}\label{Main Results}

Let $\M$ be a regular map of order $2^n$ with type $\{2^s, 2^t\}$ and $s\leq t$. For $n\leq 11$, all such regular maps $\M$ are listed in~\cite{atles}. Let $n \geq 12$. It is easy to see that if $s=1$
then $\Aut(\M) \cong C_2 \times D_{2^{n-1}}$ or $D_{2^n}$, implying $t=n-2$ or $n-1$, respectively. If $t=n-1$ then $\Aut(\M) \cong D_{2^n}$ and  $s$ is $2$ or $n-1$. Let $2\leq s,t\leq n-2$. In this section we prove that if either $s+t\leq n$, or $s+t>n$ and $s=t$, then there exists a regular map of order $2^n$ with type $\{2^s, 2^t\}$ in Theorem~\ref{existmaintheorem}, and if $s+t>n$,   $s \neq t$ and $t=n-2$ or $n-3$ then there is no regular map of order $2^n$ with type $\{2^s, 2^t\}$ in Theorem~\ref{nonexists}.

We first prove a lemma which will be used frequently in the paper.

\begin{lem}\label{relations}
Let $G=\lg \r_0, \r_1, \r_2 \rg$ such that $\r_0^2=\r_1^2=\r_2^2=(\r_0\r_2)^2=1$. Then
\begin{itemize}
\item [(1)] $[(\r_0\r_1)^2, \r_2]=[\r_0, (\r_1\r_2)^2]^{\r_2\r_1}$ and $(\r_0\r_2\r_1)^2=((\r_2\r_1)^2(\r_1\r_0)^2)^{\r_0}$;
\item [(2)] If $[(\r_0\r_1)^2, \r_2]=1$, then $[\r_0, (\r_1\r_2)^2]=[(\r_0\r_1)^4, \r_2]=[\r_0, (\r_1\r_2)^4]=1$.
 In particular, $\lg (\r_0\r_1)^2\rg\lhd G$, $\lg (\r_1\r_2)^2\rg\lhd G$, $[(\r_0\r_1)^2, (\r_1\r_2)^2]=1$ and $(\r_0\r_2\r_1)^{2^i}=(\r_2\r_1)^{2^i}(\r_1\r_0)^{2^i}$ for any positive integer $i$.
\end{itemize}

\end{lem}

\demo Since $[\r_0,\r_2]=1$, we have $[(\r_0\r_1)^2, \r_2]=[\r_0, \r_2]^{\r_1\r_0\r_1}[\r_1\r_0\r_1, \r_2]=[\r_1\r_0\r_1, \r_2]=[\r_0, \r_1\r_2\r_1]^{\r_1}=[\r_0, \r_1\r_2\r_1\r_2\r_2]^{\r_1}=[\r_0, \r_2]^{\r_1}[\r_0, (\r_1\r_2)^2]^{\r_2\r_1}=[\r_0, (\r_1\r_2)^2]^{\r_2\r_1}$ by Proposition~\ref{commutator}, and  $(\r_0\r_2\r_1)^2=\r_0\r_2\r_1\r_2\r_0\r_1=\r_0(\r_2\r_1)^2\r_1\r_0\r_1=((\r_2\r_1)^2(\r_1\r_0)^2)^{\r_0}$.

To prove (2), let $[(\r_0\r_1)^2, \r_2]=1$. By (1), $[\r_0, (\r_1\r_2)^2]=1$.
Since the relation $[(\r_0\r_1)^2, \r_2]=1$ implies that $(\r_0\r_1)^2$ and $\r_2$ commute each other, we have
$[(\r_0\r_1)^4, \r_2]=1$, and similarly, $[\r_0, (\r_1\r_2)^4]=1$.

Since $(\r_0\r_1)^{\r_1}=\r_1\r_0=(\r_0\r_1)^{-1}$ and $(\r_0\r_1)^{\r_0}=\r_1\r_0=(\r_0\r_1)^{-1}$, we have $\lg (\r_0\r_1)^2\rg\lhd G$ as $[(\r_0\r_1)^2, \r_2]=1$, and since $[\r_0,(\r_1\r_2)^2]=1$, we have $\lg (\r_1\r_2)^2\rg\lhd G$. Furthermore, $((\r_0\r_1)^2)^{(\r_1\r_2)^2}=(\r_0\r_1)^2$, that is, $[(\r_0\r_1)^2, (\r_1\r_2)^2]=1$. By (1),  $(\r_0\r_2\r_1)^2=((\r_2\r_1)^2(\r_1\r_0)^2)^{\r_0}=(\r_2\r_1)^2(\r_1\r_0)^2$, and it follows $(\r_0\r_2\r_1)^{2^i}=(\r_2\r_1)^{2^i}(\r_1\r_0)^{2^i}$ for any positive integer $i$. \hfill\qed

\begin{theorem}\label{existmaintheorem}
Let $n \geq 12$ and $2\leq s, t \leq n-2$ such that either $s+t\leq n$ or $s+t>n$ and $s=t$. Then there exists a regular map $\M$ of order $2^n$ and $\{2^s, 2^t\}$ with $\Aut(\M)=G$:

\begin{small}
\begin{itemize}
\item[(1)] for $s+t\leq n-1$, $G=\left\{
\begin{array}{ll}
\lg \r_0, \r_1, \r_2 \ |\ R_1, [(\r_0\r_1)^2, \r_2]^{2^{\frac{n-s-t-1}{2}}}\rg, & n-s-t\mbox{ odd }\\
\lg \r_0, \r_1, \r_2 \ |\ R_1, [(\r_0\r_1)^2, (\r_1\r_2)^2]^{2^{\frac{n-s-t-2}{2}}}\rg, & n-s-t\mbox{ even, }
\end{array}
\right.$

\item [(2)] for $s+t=n$, $G=\lg \r_0, \r_1, \r_2\ | \ R_2, [(\r_0\r_1)^2, \r_2], (\r_0\r_1)^{2^{s-1}}\cdot (\r_1\r_2)^{2^{t-1}}\rg$,

\item [(3)] for $s+t>n$ and $s=t$, $G=\lg \r_0, \r_1, \r_2\ | \ R_2, (\r_0\r_1)^{2^{n-t-1}}\cdot (\r_1\r_2)^{2^{n-t-1}}, [(\r_0\r_1)^2, \r_2](\r_2\r_1)^{4}\rg$,

\end{itemize}
\end{small}
where $R_1=\{\r_0^2, \r_1^2, \r_2^2, (\r_0\r_1)^{2^{s}}, (\r_1\r_2)^{2^{t}}, (\r_0\r_2)^2, [(\r_0\r_1)^4, \r_2], [\r_0,(\r_1\r_2)^4]\}$ and
$R_2=\{\r_0^2$, $\r_1^2$, $\r_2^2$,
$(\r_0\r_2)^2$, $(\r_0\r_1)^{2^{s}}, (\r_1\r_2)^{2^{t}}\}$.

\end{theorem}

\demo The part (1) is true by Proposition~\ref{$2^s,2^t$}. Now we consider parts (2) and (3), that is, the cases for $s+t=n$ and $s+t>n$ with $s=t$. To finish the proof, it suffices to show that $|G|=2^n$, $o(\r_0\r_1)=2^s$ and $o(\r_1\r_2)=2^t$.

\medskip
\f{\bf Case 1:} $s+t=n$.

In this case, we have $G=\lg \r_0, \r_1, \r_2\ | \  R_2, [(\r_0\r_1)^2, \r_2], (\r_0\r_1)^{2^{s-1}}(\r_1\r_2)^{2^{t-1}}\rg$, where $R_2=\{\r_0^2, \r_1^2, \r_2^2,(\r_0\r_2)^2, (\r_0\r_1)^{2^{s}}, (\r_1\r_2)^{2^{t}}\}$. Set
$$H=\lg  \r_0, \r_1, \r_2\ | \ \r_0^2, \r_1^2, \r_2^2, (\r_0\r_1)^{2^s},(\r_1\r_2)^{2^{t}}, (\r_0\r_2)^2, [(\r_0\r_1)^4, \r_2], [\r_0, (\r_1\r_2)^4], [(\r_0\r_1)^2, \r_2] \rg. $$

Note that $s+t<n+1$ and $(n+1)-s-t=1$. Replacing $n$ by $n+1$ in part (1), we have $|H|=2^{n+1}$, $o(\r_0\r_1)=2^s$ and  $o(\r_1\r_2)=2^{t}$ in $H$. By Lemma~\ref{relations}(2), the relations $[(\r_0\r_1)^4, \r_2]$ and $[\r_0, (\r_1\r_2)^4]$ in $H$ follow from other relations in $H$, and thus $$H=\lg  \r_0, \r_1, \r_2\ | \ \r_0^2, \r_1^2, \r_2^2, (\r_0\r_1)^{2^s},(\r_1\r_2)^{2^{t}}, (\r_0\r_2)^2, [(\r_0\r_1)^2, \r_2] \rg. $$

By Lemma~\ref{relations}(2), $\lg(\r_0\r_1)^2\rg\lhd H$ and $\lg(\r_1\r_2)^2\rg\lhd H$, and since
$o((\r_0\r_1)^{2^{s-1}})=2$ and $o((\r_1\r_2)^{2^{t-1}})=2$, we have $(\r_0\r_1)^{2^{s-1}}, (\r_1\r_2)^{2^{t-1}} \in Z(H)$, the center of $H$.

Suppose $(\r_0\r_1)^{2^{s-1}}=(\r_1\r_2)^{2^{t-1}}$. Then $H/\lg (\r_0\r_1)^{2^{s-1}}\rg\cong H_1$ with $|H_1|=2^n$, where
$$H_1=\lg  \r_0, \r_1, \r_2\ | \ \r_0^2, \r_1^2, \r_2^2, (\r_0\r_1)^{2^{s-1}},(\r_1\r_2)^{2^{t-1}}, (\r_0\r_2)^2, [(\r_0\r_1)^2, \r_2] \rg. $$
By Lemma~\ref{relations}(2), $[(\r_0\r_1)^4, \r_2]=1$ and $[\r_0, (\r_1\r_2)^4]=1$ in $H_1$, and since $(s-1)+(t-1)<n-1$ and $(n-1)-(s-1)-(t-1)=1$, part~(1) implies  $|H_1|=2^{n-1}$, a contradiction.

Thus $(\r_0\r_1)^{2^{s-1}}\not=(\r_1\r_2)^{2^{t-1}}$, that is, $o((\r_0\r_1)^{2^{s-1}}\cdot (\r_1\r_2)^{2^{t-1}})=2$ in $H$. Let $K=\lg (\r_0\r_1)^{2^{s-1}}\cdot (\r_1\r_2)^{2^{t-1}}\rg$. Then $o(\r_0\r_1K)=2^s$ and $o(\r_1\r_2K)=2^t$ in $H/K$. Clearly, $H/K\cong G$, and hence $|G|=2^n$, $o(\r_0\r_1)=2^s$ and $o(\r_1\r_2)=2^t$, as required.

\medskip
\f{\bf Case 2:} $s+t>n$ and $s=t$.

In this case, $G=\lg \r_0, \r_1, \r_2\ | \ R_2, (\r_0\r_1)^{2^{n-t-1}}(\r_1\r_2)^{2^{n-t-1}}, [(\r_0\r_1)^2, \r_2](\r_2\r_1)^{4}\rg$, where $R_2=\{\r_0^2, \r_1^2, \r_2^2,(\r_0\r_2)^2, (\r_0\r_1)^{2^{t}}, (\r_1\r_2)^{2^{t}}\}$. Set
\begin{small}$$ H=\lg  \r_0, \r_1, \r_2\ | \ \r_0^2, \r_1^2, \r_2^2, (\r_0\r_1)^{2^{n-t-1}},(\r_1\r_2)^{2^{t}}, (\r_0\r_2)^2, [(\r_0\r_1)^4, \r_2], [\r_0, (\r_1\r_2)^4], [(\r_0\r_1)^2, \r_2] \rg. $$\end{small}
Note that $(n-t-1)+t<n$ and  $n-(n-t-1)-t=1$. By part (1), we have $|H|=2^n$, $o(\r_0\r_1)=2^{n-t-1}$ and  $o(\r_1\r_2)=2^{t}$ in $H$. By Lemma~\ref{relations}(2), the relations $[(\r_0\r_1)^4, \r_2]$ and $[\r_0, (\r_1\r_2)^4]$ follows from other relations in $H$, and thus
$$H=\lg  \r_0, \r_1, \r_2\ | \ \r_0^2, \r_1^2, \r_2^2, (\r_0\r_1)^{2^{n-t-1}},(\r_1\r_2)^{2^{t}}, (\r_0\r_2)^2, [(\r_0\r_1)^2, \r_2] \rg. $$

Note that $(\r_1\r_0)^{2^{n-t-1}}=1$ and $(\r_2\r_1)^{2^{t}}=1$ in $H$. Since $s+t>n$ and $s=t$, we have $1\leq n-t-1<t$, and by Lemma~\ref{relations}(2),  $(\r_0\r_2\r_1)^{2^t}=(\r_2\r_1)^{2^t}(\r_1\r_0)^{2^t}=1$ and  $(\r_0\r_2\r_1)^{2^{n-t-1}}=(\r_2\r_1)^{2^{n-t-1}}(\r_1\r_0)^{2^{n-t-1}}=
(\r_2\r_1)^{2^{n-t-1}}$,  which implies $o(\r_0\r_2\r_1)=o(\r_2\r_1)=2^t$ in $H$ and $(\r_0\r_2\r_1)^{2^{n-t-1}}\cdot (\r_1\r_2)^{2^{n-t-1}}=1$.
Moreover, $[(\r_0\r_2\r_1)^2, \r_2]=[(\r_2\r_1)^2(\r_1\r_0)^2, \r_2]=[(\r_2\r_1)^2, \r_2]^{(\r_1\r_0)^2}[(\r_1\r_0)^2, \r_2]=((\r_1\r_2)^4)^{(\r_1\r_0)^2}=(\r_1\r_2)^4$, that is, $[(\r_0\r_2\r_1)^2, \r_2]\cdot (\r_2\r_1)^4=1$. Clearly, $H=\lg \r_0\r_2, \r_1, \r_2\rg$. Now it is easy to see that the generators $\r_0\r_2, \r_1, \r_2$ in $H$ satisfy the same relations  as $\r_0, \r_1, \r_2$ do in $G$. Thus, there is an epimorphism $\phi: G \mapsto H$ such that $\r_0^{\phi}=\r_0\r_2$, $\r_1^{\phi}=\r_1$ and $\r_2^{\phi}=\r_2$.

On the other hand, since $[(\r_0\r_1)^2, \r_2](\r_2\r_1)^{4}=1$ in $G$,  Lemma~\ref{relations}(1) implies $(\r_1\r_2)^4=[(\r_0\r_1)^2, \r_2]=[\r_0, (\r_1\r_2)^2]^{\r_2\r_1}$, that is, $[\r_0, (\r_1\r_2)^2]=(\r_1\r_2)^4$ in $G$.
Thus, $[\r_0\r_1, (\r_1\r_2)^2]=[\r_0, (\r_1\r_2)^2]^{\r_1}[\r_1, (\r_1\r_2)^2]=[(\r_1\r_2)^4]^{\r_1}(\r_1\r_2)^4=1$ and $[(\r_0\r_1)^2, (\r_1\r_2)^2]=1$.
Note that $(\r_0\r_2\r_1)^2=\r_2(\r_0\r_1)^2\r_2(\r_2\r_1)^2=(\r_0\r_1)^2[(\r_0\r_1)^2,\r_2](\r_2\r_1)^2=(\r_0\r_1)^2(\r_1\r_2)^4(\r_2\r_1)^2=
(\r_0\r_1)^2(\r_1\r_2)^2$. This implies $(\r_0\r_2\r_1)^{2^{n-t-1}}=1$ because  $(\r_0\r_1)^{2^{n-t-1}}\cdot (\r_1\r_2)^{2^{n-t-1}}=1$ in $G$ and $n-t-1\geq 1$.
Furthermore, $[(\r_0\r_2\r_1)^2, \r_2]=[(\r_0\r_1)^2(\r_1\r_2)^2, \r_2]=[(\r_0\r_1)^2, \r_2]^{(\r_1\r_2)^2}$
$[(\r_1\r_2)^2, \r_2]=(\r_1\r_2)^4(\r_2\r_1)^4=1$. Thus, the generators $\r_0\r_2, \r_1, \r_2$ in $G$ satisfy the same relations as $\r_0, \r_1, \r_2$ do  in $H$, which implies that  there is an epimorphism $\varphi: H \mapsto G$ such that $\r_0^{\varphi}=\r_0\r_2$, $\r_1^{\varphi}=\r_1$ and $\r_2^{\varphi}=\r_2$.

Since $|H|=2^n$, both $\phi$ and $\varphi$ are isomorphisms, and hence $|G|=2^n$. Since   $o(\r_1\r_2)=2^{t}$ in $H$ and $(\r_1\r_2)^{\varphi}=\r_1\r_2$, we have $o(\r_1\r_2)=2^{t}$ in $G$, and since $o(\r_0\r_2\r_1)=2^t$ in $H$ and $(\r_0\r_2\r_1)^{\varphi}=\r_0\r_1$, we have $o(\r_0\r_1)=2^t$ in $G$, as required.  \hfill\qed

\begin{theorem}\label{nonexists}
Let $n \geq 12$ and $2\leq s,t \leq n-2$ such that $s+t>n$ and $s<t$. Then there is no regular map of order $2^n$ with type $\{2^s, 2^t\}$ for $t=n-2$ or $n-3$.
\end{theorem}

\demo Suppose to the contrary that $\M$ is a regular map of order $2^n$ with type $\{2^s, 2^t\}$ for $t=n-2$ or $n-3$. Let $G=\Aut(\M)$. Then $|G|=2^n$ and $\Aut(\M)$ is generated by three involutions, say $\r_0, \r_1$ and $\r_2$, such that $\r_0\r_2=\r_2\r_0$, $o(\r_0\r_1)=2^s$ and $o(\r_1\r_2)=2^t$. Since $|G|=2^n$, $G$ cannot be generated by any two of $\r_0$, $\r_1$ and $\r_2$, that is, $\{\r_0, \r_1,\r_2\}$ is a minimal generating set of $G$. In particular, $|\lg \r_1, \r_2\rg|=2^{t+1}$ and $|\lg \r_0, \r_1\rg|=2^{s+1}$. Note that    $[\r_0,\r_1]=[\r_0,\r_2\r_2\r_1]=[\r_0,\r_2\r_1]=(\r_1\r_2)^{\r_0}\r_2\r_1$, that is, $(\r_1\r_2)^{\r_0}=(\r_0\r_1)^2\r_1\r_2$.

Let $t=n-2$. Then $o(\r_1\r_2)=2^{n-2}$ and $|\lg \r_1, \r_2\rg|=2^{n-1}$. Since  $|G|=2^n$, we have $\lg \r_1, \r_2\rg \unlhd G$, and hence $\lg \r_1\r_2\rg \unlhd G$, because $\lg \r_1\r_2\rg$ is characteristic in $\lg \r_1, \r_2\rg$. It follows $(\r_1\r_2)^{\r_0}=(\r_1\r_2)^i$ for some $i\in\mz_{2^{n-2}}$. Since $o(\r_0)=2$ and $o(\r_1\r_2)=2^{n-2}$, we have $i^2=1$ in $\mz_{2^{n-2}}^*$, where $\mz_{2^{n-2}}^*$ is the multiplicative group of $\mz_{2^{n-2}}$ consisting of numbers of $\mz_{2^{n-2}}$ coprime to $2$.
It is well known that $\mz_{2^{n-2}}^*\cong \mz_{2^{n-4}}\times \mz_2$, and so the equation $i^2=1$ has four solutions in $\mz_{2^{n-2}}^*$, that is, $i=1$, $-1$, $2^{n-3}+1$ or $2^{n-3}-1$. It follows $(\r_1\r_2)^{\r_0}=\r_1\r_2, \r_2\r_1, \r_1\r_2(\r_1\r_2)^{2^{n-3}}$ or $\r_2\r_1(\r_1\r_2)^{2^{n-3}}$, and hence  $(\r_0\r_1)^2=[\r_0,\r_1]=(\r_1\r_2)^{\r_0}\r_2\r_1=1,(\r_2\r_1)^2, (\r_1\r_2)^{2^{n-3}}$ or $(\r_2\r_1)^2(\r_1\r_2)^{2^{n-3}}$, which implies $s=1,2$ or $s=t=2^{n-2}$, but this is impossible because $s+t>n$ and $s<t$.

Let $t=n-3$. Then $o(\r_1\r_2)=2^{n-3}$ and $|\lg \r_1, \r_2\rg|=2^{n-2}$. Let $A=\lg \r_0\r_1, \r_1\r_2\rg$ and $B=\lg (\r_0\r_1)^2, \r_1\r_2\rg$. Then $A \unlhd G$ and $|G:A|\leq 2$. Since $\{\r_0,\r_1,\r_2\}$ is a minimal generating set, by Proposition~\ref{burnside}, $G$ has rank $3$, that is, $d(G)=3$. This implies that $d(A)=2$ and $|A|=2^{n-1}$.
Clearly, $(\r_1\r_2)^{\r_0\r_1}=(\r_1\r_0)^2(\r_2\r_1) \in B$, and hence $B \unlhd A$. Thus $|A/B| \leq 2$. By Proposition~\ref{burnside}, $(\r_0\r_1)^2 \in \mho_{1}(A) \leq \Phi(A)$.
If $B=A$, then $A/\Phi(A) =B/\Phi(A)=\lg \r_1\r_2\Phi(A)\rg$, and by Proposition~\ref{frattini}, $A=\lg \r_1\r_2,\Phi(A)\rg=\lg \r_1\r_2\rg$, contradicting  $d(A)=2$. Thus $B\not=A$ and $|B|=2^{n-2}$.

Recall $(\r_1\r_2)^{\r_0}=(\r_0\r_1)^2\r_1\r_2\in B$.
Since $|B|=2^{n-2}$ and $o(\r_1\r_2)=o((\r_1\r_2)^{\r_0})=2^{n-3}$, we have $\lg (\r_1\r_2)^2\rg =\lg ((\r_1\r_2)^{\r_0})^2 \rg$, and since $o((\r_1\r_2)^2)=2^{n-4}$, we have  $((\r_1\r_2)^{\r_0})^2=(\r_1\r_2)^{2i}$ for some $i \in \mz_{2^{n-4}}$. Now $(\r_1\r_2)^2=((\r_1\r_2)^2)^{\r_0^2}=((\r_1\r_2)^{2i})^{\r_0}=(((\r_1\r_2)^{\r_0})^2)^{i}=((\r_1\r_2)^2)^{i^2}$,
which implies
$i^2-1=0$ in $\mz_{2^{n-4}}^*$. It follows $i=\pm1$ or $2^{n-5}\pm1$.

Since $\r_1\r_2$ and $(\r_1\r_2)^{\r_0}$ are elements of order $2^{n-3}$ in $B$ and $|B|=2^{n-2}$, Corollary~\ref{maximalc} implies $(\r_1\r_2)^{(\r_1\r_2)^{\r_0}}=(\r_1\r_2)$ or $(\r_1\r_2)^{1+2^{n-4}}$, and since  $(\r_0\r_1)^2=(\r_1\r_2)^{\r_0}\r_2\r_1$, we have $(\r_1\r_2)^{(\r_0\r_1)^2}=(\r_1\r_2)$ or $(\r_1\r_2)^{1+2^{n-4}}$. It follows
$((\r_1\r_2)^{\r_0})^2=((\r_0\r_1)^2\r_1\r_2)((\r_0\r_1)^2\r_1\r_2)=
(\r_0\r_1)^4(\r_1\r_2)^{(\r_0\r_1)^2}\r_1\r_2=(\r_0\r_1)^4(\r_1\r_2)^2$ or $(\r_0\r_1)^4(\r_1\r_2)^{2+2^{n-4}}$.

Let $((\r_1\r_2)^{\r_0})^2=(\r_0\r_1)^4(\r_1\r_2)^2$. Then $(\r_1\r_2)^{2i}=(\r_0\r_1)^4(\r_1\r_2)^2$, and since $i=\pm1$ or $2^{n-5}\pm1$, we have $(\r_0\r_1)^4=(\r_1\r_2)^{2i-2}=1$, $(\r_2\r_1)^4$, $(\r_1\r_2)^{2^{n-4}}$ or $(\r_2\r_1)^4(\r_1\r_2)^{2^{n-4}}$. It follows $s=2,3$ or $s=t=2^{n-3}$, each of which is impossible because $s+t>n$ and $s<t$.

Let $((\r_1\r_2)^{\r_0})^2=(\r_0\r_1)^4(\r_1\r_2)^{2+2^{n-4}}$.
Since $o(\r_1\r_2)=2^{n-3}$, we have  $(\r_0\r_1)^4=(\r_1\r_2)^{2^{n-4}+2i-2}=(\r_1\r_2)^{2^{n-4}}$, $(\r_2\r_1)^4(\r_1\r_2)^{2^{n-4}}$, $1$ or $(\r_2\r_1)^4$. This means that $s=2,3$ or $s=t=2^{n-3}$, each of which is impossible because $s+t>n$ and $s<t$.
\hfill\qed

Based on Theorem~\ref{nonexists}, we would like to propose the following conjecture.

\begin{conjecture}
For any positive integers $n, s, t$ such that $n \geq 12$, $2 \leq s, t \leq n-2$, $s+t>n$ and $s \ne t$, there is no regular map of order $2^n$ with type $\{2^s, 2^t\}$.
\end{conjecture}

A computation with {\sc Magma} shows that the conjecture is true for $n\leq 12$.

\section{Regular maps of order $2^n$ with certain types}\label{Main Results-classifications}

In this section, we classify regular maps of order $2^n$ with types $\{2^{n-2},2^{n-2}\}$
and $\{2^{n-3},2^{n-3}\}$ in Theorem~\ref{maintheorem2}. It appears that these regular maps have the following automorphism groups $G_i$ for $1\leq i\leq 6$. Let $n \geq 12$ and let $$R_1=\{ \r_0^2, \r_1^2, \r_2^2, (\r_0\r_2)^2\} \mbox{ and } R_2=\{ \r_0^2, \r_1^2, \r_2^2, (\r_0\r_2)^2,(\r_0\r_1)^{2^{n-3}}, (\r_1\r_2)^{2^{n-3}}\}.$$ We define six groups as following:

\begin{itemize}\setlength{\parskip}{-3pt}
  \item [$G_1$]$=\lg \r_0, \r_1, \r_2 \ |\ R_1, (\r_0\r_1)^{2^{n-2}}, (\r_1\r_2)^{2^{n-2}}, (\r_0\r_1)^2(\r_1\r_2)^2 \rg$,
  \item [$G_2$]$=\lg \r_0, \r_1, \r_2 \ |\ R_1, (\r_0\r_1)^{2^{n-2}}, (\r_1\r_2)^{2^{n-2}}, (\r_0\r_1)^2(\r_1\r_2)^2(\r_1\r_2)^{2^{n-3}}\rg$,
  \item [$G_3$]$=\lg \r_0, \r_1, \r_2 \ |\ R_2, (\r_0\r_1)^4(\r_1\r_2)^4, [(\r_0\r_1)^2,\r_2](\r_0\r_1)^4\rg$,
  \item [$G_4$]$=\lg \r_0, \r_1, \r_2 \ |\ R_2, (\r_0\r_1)^4(\r_1\r_2)^4, [(\r_0\r_1)^2,\r_2](\r_0\r_1)^4(\r_1\r_2)^{2^{n-4}}\rg,$
  \item [$G_5$]$=\lg \r_0, \r_1, \r_2 \ |\ R_2, (\r_0\r_1)^4(\r_1\r_2)^4(\r_1\r_2)^{2^{n-4}}, [(\r_0\r_1)^2,\r_2](\r_0\r_1)^4(\r_1\r_2)^{2^{n-4}}\rg$,
  \item [$G_6$]$=\lg \r_0, \r_1, \r_2 \ |\ R_2, (\r_0\r_1)^4(\r_1\r_2)^4(\r_1\r_2)^{2^{n-4}}, [(\r_0\r_1)^2,\r_2](\r_0\r_1)^4\rg$.
 \end{itemize}

We first prove that the groups $G_i$ for $1\leq i\leq 5$ have order $2^n$. This is also true for $G_6$, but the proof is quite different with that given in Theorem~\ref{maintheorem2}.

\begin{lem}\label{grouporder}
The groups $G_i$ for $1\leq i\leq 5$ have order $2^n$, and
the listed exponents are the true orders of the corresponding elements in each $G_i$. Furthermore, $|G_6|\leq 2^n$.
\end{lem}

\demo Note that $(n-2)+(n-2)>n$ and $(n-3)+(n-3)>n$. By taking $(s, t, n)=(n-2, n-2, n),(n-3,n-3, n)$ in Theorem~\ref{existmaintheorem}(3), we know that both groups

\begin{small}
\begin{itemize}\setlength{\parskip}{0pt}
  \item [$H_1$]$=\lg \r_0, \r_1, \r_2 \ |\ \r_0^2, \r_1^2, \r_2^2, (\r_0\r_1)^{2^{n-2}}, (\r_1\r_2)^{2^{n-2}}, (\r_0\r_2)^2, (\r_0\r_1)^{2}(\r_1\r_2)^{2}, [(\r_0\r_1)^2,\r_2](\r_2\r_1)^4 \rg$,
   \item [$H_3$]$=\lg \r_0, \r_1, \r_2 \ |\ \r_0^2, \r_1^2, \r_2^2, (\r_0\r_1)^{2^{n-3}}, (\r_1\r_2)^{2^{n-3}}, (\r_0\r_2)^2, (\r_0\r_1)^4(\r_1\r_2)^4, [(\r_0\r_1)^2,\r_2](\r_2\r_1)^4\rg$,
  \end{itemize}
\end{small}
have order $2^n$ and
the listed exponents are the true orders of the corresponding elements in $H_1$ and $H_3$. For $H_1$, we have $(\r_0\r_1)^2=(\r_1\r_2)^{-2}$, and since $[(\r_0\r_1)^2, \r_2]=(\r_0\r_1)^{-2}\r_2(\r_0\r_1)^2\r_2=(\r_1\r_2)^2\r_2(\r_1\r_2)^{-2}\r_2=(\r_1\r_2)^4$, the relation $[(\r_0\r_1)^2, \r_2](\r_2\r_1)^4$ in $H_1$ is redundant and so $G_1=H_1$.
For $H_3$, we have $(\r_0\r_1)^4=(\r_2\r_1)^{4}$, and hence $[(\r_0\r_1)^2, \r_2](\r_2\r_1)^4=1$ if and only if $[(\r_0\r_1)^2, \r_2](\r_0\r_1)^4=1$. It follows  $G_3=H_3$.

For $G_2$, we have $(\r_0\r_1)^2=(\r_2\r_1)^{2+2^{n-3}}$ and  $((\r_0\r_1)^2)^{\r_2}=(\r_1\r_2)^{2+2^{n-3}}$. By Proposition~\ref{derived},
$G_2'=\lg [\r_0, \r_1], [\r_1, \r_2], [\r_0,\r_1]^{\r_2}\rg=\lg (\r_0\r_1)^2, (\r_1\r_2)^2, ((\r_0\r_1)^2)^{\r_2}\rg=\lg (\r_1\r_2)^2\rg$, and so $|G_2'| = o((\r_1\r_2)^2) \leq 2^{n-3}$.
Note that $G_2/G_2'$ is abelian and is generated by three involutions. Thus $|G_2/G_2'| \leq 2^3$ and hence $|G_2|=|G_2/G_2'| \cdot |G_2'|\leq 2^n$.

By taking $(s, t, n)=(2, n-2, n)$ in Theorem~\ref{existmaintheorem}(2), we know that the group
$$H_2=\lg \r_0, \r_1, \r_2 \ | \ \r_0^2, \r_1^2, \r_2^2, (\r_0\r_1)^4, (\r_1\r_2)^{2^{n-2}}, (\r_0\r_2)^2, (\r_0\r_1)^2(\r_1\r_2)^{2^{n-3}},
[(\r_0\r_1)^2, \r_2]\rg$$
has order $2^n$ and
the listed exponents are the true orders of the corresponding elements in $H_2$. Since $[(\r_0\r_1)^2, \r_2]=1$ in $H_2$, Lemma~\ref{relations}(2) implies that
$(\r_0\r_2\r_1)^2=(\r_2\r_1)^2(\r_1\r_0)^2=(\r_2\r_1)^2(\r_2\r_1)^{2^{n-3}}$ in $H_2$. It follows $o(\r_0\r_2\r_1)=2^{n-2}$ and $(\r_0\r_2\r_1)^2(\r_1\r_2)^2(\r_1\r_2)^{2^{n-3}}=1$ in $H_2$. So the generators $\r_0\r_2, \r_1, \r_2$ in $H_2$ satisfy the same relations  as $\r_0, \r_1, \r_2$ do in $G_2$, and hence there is an epimorphism $\phi: G_2 \mapsto H_2$ such that $\r_0^{\phi}=\r_0\r_2$, $\r_1^{\phi}=\r_1$ and $\r_2^{\phi}=\r_2$. Since $|H_2|=2^n$ and $|G_2| \leq 2^n$, we have $|G_2|=2^n$, and since $(\r_0\r_1)^\phi=\r_0\r_2\r_1$, we have $o(\r_0\r_1)=2^{n-2}$ in $G_2$.

Let $G=G_4, G_5$ or $G_6$. Then $(\r_0\r_1)^4=(\r_2\r_1)^4$ or $(\r_2\r_1)^4(\r_2\r_1)^{2^{n-4}}$, and since $(\r_1\r_2)^{2^{n-3}}=1$, we have $(\r_0\r_1)^{2^{n-4}}=((\r_0\r_1)^4)^{2^{n-8}}=(\r_2\r_1)^{2^{n-4}}$. Let $K=\lg (\r_1\r_2)^{2^{n-4}}\rg$. Then $|K|\leq 2$ and $K\lhd G$, because $(\r_0\r_1)^{\r_0}=(\r_0\r_1)^{-1}$, $(\r_0\r_1)^{\r_1}=(\r_0\r_1)^{-1}$ and $(\r_1\r_2)^{\r_2}=(\r_1\r_2)^{-1}$. The three generators $\r_0K, \r_1K, \r_2K$ in $G/K$ satisfy the same relations as $\r_0, \r_1, \r_2$ in $G_3$, where $n$ is replaced by $n-1$,  and hence $|G/K| \leq 2^{n-1}$ (here, we need check $|G_3|=2^{11}$ for $n=11$ and this can be done by {\sc Magma}). It follows $|G| \leq 2^{n-1}\cdot 2=2^n$.

Set $H_4=\lg \r_0, \r_1, \r_2 \ | \ \r_0^2, \r_1^2, \r_2^2, (\r_0\r_1)^4, (\r_1\r_2)^{2^{n-3}}, (\r_0\r_2)^2, [(\r_0\r_1)^2, \r_2](\r_2\r_1)^{2^{n-4}}\rg$. Note that $(\r_2\r_1)^{2^{n-4}}=(\r_1\r_2)^{2^{n-4}}$ in $H_4$. Now $H_4=L$ with $L$ given in Proposition~\ref{$2^s,2^t$}, and hence $|H_4|=2^n$ and  the listed exponents are the true orders of the corresponding elements in $H_4$.
Since $[(\r_0\r_1)^2, \r_2]=(\r_1\r_2)^{2^{n-4}}$ and $(\r_0\r_1)^4=1$ in $H_4$,  Proposition~\ref{commutator} implies $[(\r_0\r_1)^2, \r_1\r_2]=[(\r_0\r_1)^2, \r_2][(\r_0\r_1)^2, \r_1]^{\r_2}=(\r_1\r_2)^{2^{n-4}}$ and hence $[(\r_0\r_1)^2, (\r_1\r_2)^2]=[(\r_0\r_1)^2, \r_1\r_2][(\r_0\r_1)^2, \r_1\r_2]^{\r_1\r_2}=(\r_1\r_2)^{2^{n-3}}=1$ in $H_4$. By Lemma~\ref{relations}(1), $[\r_0, (\r_1\r_2)^2]=[(\r_0\r_1)^2, \r_2]^{\r_1\r_2}=(\r_1\r_2)^{2^{n-4}}$ and hence $[\r_0, (\r_1\r_2)^4]
=[\r_0, (\r_1\r_2)^2][\r_0, (\r_1\r_2)^2]^{(\r_1\r_2)^2}=(\r_1\r_2)^{2^{n-4}}(\r_1\r_2)^{2^{n-4}}=1$.
Again by Lemma~\ref{relations}(1), $(\r_0\r_2\r_1)^4=((\r_2\r_1)^4(\r_1\r_0)^4)^{\r_0}=((\r_2\r_1)^4)^{\r_0}=(\r_2\r_1)^4$ in $H_4$, which implies $o(\r_0\r_2\r_1)=2^{n-3}$ and $(\r_0\r_2\r_1)^4(\r_1\r_2)^4=1$.

By Lemma~\ref{relations}(1), $(\r_0\r_2\r_1)^2=((\r_2\r_1)^2(\r_1\r_0)^2)^{\r_0}$,
and since $[(\r_2\r_1)^2(\r_1\r_0)^2, \r_2])=[(\r_2\r_1)^2, \r_2]^{(\r_1\r_0)^2}[(\r_1\r_0)^2, \r_2]
=((\r_1\r_2)^4)^{(\r_1\r_0)^2}(\r_1\r_2)^{2^{n-4}}=(\r_1\r_2)^4(\r_1\r_2)^{2^{n-4}}$, one may have $[(\r_0\r_2\r_1)^2, \r_2]=
[(\r_2\r_1)^2(\r_1\r_0)^2, \r_2])^{\r_0}=(\r_1\r_2)^4(\r_1\r_2)^{2^{n-4}}$.
Since $(\r_0\r_2\r_1)^4(\r_1\r_2)^4=1$, we have $[(\r_0\r_2\r_1)^2, \r_2](\r_0\r_2\r_1)^4(\r_1\r_2)^{2^{n-4}}=1$. So the generators $\r_0\r_2, \r_1, \r_2$ in $H_4$ satisfy the same relations  as $\r_0, \r_1, \r_2$ do in $G_4$, and hence there is an epimorphism $\phi: G_4 \mapsto H_4$ such that $\r_0^{\phi}=\r_0\r_2$, $\r_1^{\phi}=\r_1$ and $\r_2^{\phi}=\r_2$. Since $|H_4|=2^n$ and $|G_4| \leq 2^n$, we have $|G_4|=2^n$, and since $(\r_0\r_1)^\phi=\r_0\r_2\r_1$, we have $o(\r_0\r_1)=o(\r_1\r_2)=2^{n-3}$ in $G_4$.

By taking $(s, t, n)=(3, n-3, n)$ in Theorem~\ref{existmaintheorem}(2), we know that the group
$$H_5=\lg \r_0, \r_1, \r_2 \ | \ \r_0^2, \r_1^2, \r_2^2, (\r_0\r_1)^8, (\r_1\r_2)^{2^{n-3}}, (\r_0\r_2)^2,[(\r_0\r_1)^2, \r_2], (\r_0\r_1)^4(\r_1\r_2)^{2^{n-4}}\rg$$
has order $2^n$ and the listed exponents are the true orders of the corresponding elements in $H_5$. Since $[(\r_0\r_1)^2, \r_2]=1$ and $(\r_0\r_1)^4=(\r_1\r_0)^4$ in $H_5$, by Lemma~\ref{relations}(2) we have   $(\r_0\r_2\r_1)^4=(\r_2\r_1)^4(\r_0\r_1)^4=(\r_2\r_1)^4(\r_2\r_1)^{2^{n-4}}$ and  $(\r_0\r_2\r_1)^8=(\r_2\r_1)^8(\r_0\r_1)^8=(\r_2\r_1)^8$, implying  $o(\r_0\r_2\r_1)=o(\r_1\r_2)=2^{n-3}$ and $(\r_0\r_2\r_1)^4(\r_1\r_2)^4(\r_1\r_2)^{2^{n-4}}=1$ in $H_5$. Again by Lemma~\ref{relations}(2), $[\r_0,(\r_1\r_2)^2]=1$ and hence
$[(\r_0\r_2\r_1)^2, \r_2]=[(\r_2\r_1)^2(\r_0\r_1)^2, \r_2]=[(\r_2\r_1)^2, \r_2]^{(\r_0\r_1)^2}[(\r_0\r_1)^2, \r_2]=(\r_1\r_2)^4$ in $H_5$, which implies $[(\r_0\r_2\r_1)^2, \r_2](\r_0\r_2\r_1)^4(\r_1\r_2)^{2^{n-4}}$
$=1$.
Thus the generators $\r_0\r_2, \r_1, \r_2$ in $H_5$ satisfy the same relations as  $\r_0, \r_1, \r_2$ do in $G_5$, and so there is an epimorphism $\phi: G_5 \mapsto H_5$ such that $\r_0^{\phi}=\r_0\r_2$, $\r_1^{\phi}=\r_1$ and $\r_2^{\phi}=\r_2$. Since $|H_5|=2^n$ and $|G| \leq 2^n$, we have $|G_5|=2^n$, and since $(\r_0\r_1)^\phi=\r_0\r_2\r_1$, we have $o(\r_0\r_1)=2^{n-3}$ and $o(\r_1\r_2)=2^{n-3}$ in $G_5$.\hfill\qed

The following concerns a quotient of a regular map of order a $2$-power.

\begin{lem}\label{quotient}
Let $n \geq 12$ and $s+t>n$ with $2\leq s,t\leq n-2$. Let $G=\lg \r_0, \r_1, \r_2\rg$ such that $|G|=2^n$, $o(\r_0)=o(\r_1)=o(\r_2)=o(\r_0\r_2)=2$, $o(\r_0\r_1)=2^s$ and $o(\r_1\r_2)=2^t$. Then $N =\lg (\r_1\r_2)^{2^{t-1}} \rg\unlhd G$ and $G/N=\lg \r_0N,\r_1N,\r_2N\rg$. Furthermore, $G/N=2^{n-1}$, $o(\r_0N)=o(\r_1N)=o(\r_2N)=o(\r_0\r_2N)=2$,  $o(\r_0\r_1N)=2^{s-1}$ and $o(\r_1\r_2N)=2^{t-1}$.
\end{lem}

\demo Since $|\lg \r_0\r_1\rg \lg \r_1\r_2\rg|=\frac{|\lg \r_0\r_1\rg||\lg \r_1\r_2\rg|}{|\lg \r_0\r_1\rg\cap \lg \r_1\r_2\rg|}=\frac{2^{s+t}}{|\lg \r_0\r_1\rg\cap \lg \r_1\r_2\rg|}>\frac{2^n}{|\lg \r_0\r_1\rg\cap \lg \r_1\r_2\rg|}=\frac{|G|}{|\lg \r_0\r_1\rg\cap \lg \r_1\r_2\rg|}$, we have $\lg \r_0\r_1\rg\cap \lg \r_1\r_2\rg\not=1$. Note that $\r_0$ and $\r_1$ normalize $\lg \r_0\r_1\rg$, and $\r_1$ and $\r_2$ normalize $\lg \r_1\r_2\rg$. Since $G=\lg \r_0,\r_1,\r_2\rg$, we have $\lg \r_0\r_1\rg\cap \lg \r_1\r_2\rg\unlhd G$, and since every subgroup of $\lg \r_0\r_1\rg\cap \lg \r_1\r_2\rg$ is characteristic in the cyclic group $\lg \r_0\r_1\rg\cap \lg \r_1\r_2\rg$, it is normal in $G$. In particular, $N=\lg (\r_0\r_1)^{2^{s-1}}\rg=\lg (\r_1\r_2)^{2^{t-1}}\rg\cong \mz_2$ and $N\unlhd G$. It follows $|G/N|=2^{n-1}$, $G/N=\lg \r_0N,\r_1N,\r_2N\rg$, $o(\r_0\r_1N)=2^{s-1}$ and $o(\r_1\r_2N)=2^{t-1}$.

Note that $\lg \r_0,\r_1\rg\cong D_{2^{s+1}}$, $\lg \r_1,\r_2\rg\cong D_{2^{t+1}}$ and $\lg\r_0,\r_2\rg\cong\mz_2\times \mz_2$. Since $s,t\leq n-2$ and $|G|=2^n$, $\{\r_0,\r_1,\r_2\}$ is a minimal generating set of $G$. It follows $\r_0,\r_0\r_2,\r_2\not\in N$, and hence   $o(\r_0N)=o(\r_2N)=o(\r_0\r_2N)=2$. If $\r_1\in N$, then $\lg \r_1\rg\unlhd G$ and hence $|G|\leq |\lg \r_0,\r_2\rg||\lg \r_1\rg|\leq 2^3$, a contradiction. Thus $\r_1\not\in N$ and $o(\r_1N)=2$. \hfill\qed

\begin{theorem}\label{maintheorem2}
Let $n \geq 12$ and let $\M$ be a regular map of order $2^n$. Then
\begin{enumerate}
\item[(1)] $\M$ has type  $\{2^{n-2},2^{n-2}\}$ if and only if $\Aut(\M)\cong G_1$ or $G_2$;
\item[(2)] $\M$ has type  $\{2^{n-3},2^{n-3}\}$ if and only if $\Aut(\M)\cong G_3, G_4, G_5$ or $G_6$.
\end{enumerate}
\end{theorem}

\demo By Lemma~\ref{grouporder}, $G_1$ and $G_2$ are automorphism groups of regular maps of order $2^n$ and  type $\{2^{n-2},2^{n-2}\}$. For the necessity part in (1), let $G=\lg \r_0,\r_1,\r_2\rg$ be the automorphism group of a regular map of order $2^{n}$ with type $\{2^{n-2}, 2^{n-2}\}$.  Then $o(\r_0)=o(\r_1)=o(\r_2)=o(\r_0\r_2)=2$ and $o(\r_0\r_1)=o(\r_1\r_2)=2^{n-2}$. To finish the proof of part (1), we only need to show $G=G_1$ or $G_2$, that is, to show $(\r_0\r_1)^2(\r_1\r_2)^2=1$ or $(\r_0\r_1)^2(\r_1\r_2)^2(\r_1\r_2)^{2^{n-3}}=1$ in $G$. This is true for $n=12$ by {\sc Magma}. Let us begin by induction on $n$

Assume $n\geq 13$. Take $N=\lg (\r_1\r_2)^{2^{n-3}}\rg$. By Lemma~\ref{quotient}, $N\unlhd G$ and $G/N$ is the automorphism group of some regular map of order $2^{n-1}$ with type $\{2^{n-3}, 2^{n-3}\}$. Write $\olg=G/N$ and $\overline{x}=xN$ for any $x\in G$. Since $|\olg|=2^{n-1}$, the induction hypothesis implies that we may assume $\olg=\olg_1$ or $\olg_2$, where

\begin{itemize}
  \item [$\olg_1$]$=\lg \overline{\r_0}, \overline{\r_1}, \overline{\r_2} \ |\ \overline{\r_0}^2, \overline{\r_1}^2, \overline{\r_2}^2, (\overline{\r_0}\overline{\r_1})^{2^{n-3}}, (\overline{\r_1}\overline{\r_2})^{2^{n-3}}, (\overline{\r_0}\overline{\r_2})^2, (\overline{\r_0}\overline{\r_1})^2(\overline{\r_1}\overline{\r_2})^2\rg$,

  \item [$\olg_2$]$=\lg \overline{\r_0}, \overline{\r_1}, \overline{\r_2} \ |\ \overline{\r_0}^2, \overline{\r_1}^2, \overline{\r_2}^2, (\overline{\r_0}\overline{\r_1})^{2^{n-3}}, (\overline{\r_1}\overline{\r_2})^{2^{n-3}}, (\overline{\r_0}\overline{\r_2})^2,(\overline{\r_0}\overline{\r_1})^2(\overline{\r_1}\overline{\r_2})^2(\overline{\r_1}\overline{\r_2})^{2^{n-4}}\rg$.
\end{itemize}

Suppose $\olg=\olg_2$. Since $\lg (\r_1\r_2)^{2^{n-3}}\rg \cong \mathbb{Z}_2$, we have $(\r_0\r_1)^2(\r_1\r_2)^2(\r_1\r_2)^{2^{n-4}}=1$ or $(\r_1\r_2)^{2^{n-3}}$, that is,  $(\r_0\r_1)^2(\r_1\r_2)^2=(\r_1\r_2)^{\d \cdot2^{n-4}}$ with $\d =1$ or $-1$. It follows $(\r_0\r_1)^2=(\r_1\r_2)^{\d \cdot2^{n-4}-2} \in \lg (\r_1\r_2)^2\rg$, and hence $[(\r_0\r_1)^2, (\r_1\r_2)^2]=1$ and $((\r_0\r_1)^2)^{\r_2}=(\r_0\r_1)^{-2}$. Thus $(\r_0\r_1)^4(\r_1\r_2)^4=((\r_0\r_1)^2(\r_1\r_2)^2)^2=(\r_1\r_2)^{2^{n-3}}$ and $(\r_0\r_1)^8(\r_1\r_2)^8=1$, which implies  $(\r_1\r_2)^{\d \cdot2^{n-4}}=(\r_1\r_0)^{\d \cdot2^{n-4}}$.

Since $(\r_1\r_2)^2=(\r_1\r_0)^2(\r_1\r_2)^{\d \cdot2^{n-4}}=(\r_1\r_0)^{\d \cdot2^{n-4}+2}$, we have $((\r_2\r_1)^2)^{\r_0}=(\r_2\r_1)^{-2}$.
Note that $[(\r_0\r_1)^2, \r_2]=(\r_1\r_0)^2((\r_0\r_1)^2)^{\r_2}=(\r_1\r_0)^4$
and $[\r_0, (\r_1\r_2)^2]=((\r_2\r_1)^2)^{\r_0}(\r_1\r_2)^2$
$=(\r_1\r_2)^4$. By Lemma~\ref{relations}(1), $[\r_0, (\r_1\r_2)^2]=[(\r_0\r_1)^2, \r_2]^{\r_1\r_2}$ and so $(\r_1\r_0)^4=(\r_1\r_2)^4$. It follows  $(\r_1\r_2)^{2^{n-3}}=(\r_0\r_1)^4(\r_1\r_2)^4
=1$, contrary to $o(\r_1\r_2)=2^{n-2}$.

Now we have $\olg=\olg_1$. Since $N=\langle (\r_1\r_2)^{2^{n-3}}\rangle\cong\mz_2$, we have $(\r_0\r_1)^2(\r_1\r_2)^2=1$ or $(\r_1\r_2)^{2^{n-3}}$. For the latter,  $(\r_0\r_1)^2(\r_1\r_2)^2(\r_1\r_2)^{2^{n-3}}=1$. Thus $G\cong G_1$ or $G_2$.

\medskip

For the sufficiency part in (2), by Lemma~\ref{grouporder} we only need to show that $G_6$ has order $2^{n}$ and the listed exponents are the true orders of the corresponding elements. Let
$$H_6=\lg \r_0, \r_1, \r_2 \ | \ \r_0^2, \r_1^2, \r_2^2, (\r_0\r_1)^8, (\r_1\r_2)^{2^{n-3}}, (\r_0\r_2)^2, (\r_0\r_1)^4(\r_1\r_2)^{2^{n-4}}, [(\r_0\r_1)^2, \r_2](\r_0\r_1)^4\rg.$$

Note that  $\lg \r_0\r_1\rg\cap\lg \r_1\r_2\rg\unlhd H_6$ because $\lg \r_0\r_1\rg\cap\lg \r_1\r_2\rg$ is normalized by $\r_0,\r_1$ and $\r_2$. Since  $(\r_0\r_1)^4(\r_1\r_2)^{2^{n-4}}=1$, we have $\lg (\r_0\r_1) ^4\rg \unlhd H_6$ and $H_6/\lg (\r_0\r_1)^4\rg\cong L_6$, where
$$L_6=\lg \r_0, \r_1, \r_2 \ |\ \r_0^2, \r_1^2, \r_2^2, (\r_0\r_1)^{2^{2}}, (\r_1\r_2)^{2^{n-4}}, (\r_0\r_2)^2, [(\r_0\r_1)^2,\r_2]\rg.$$
Note that $[(\r_0\r_1)^2,\r_2]=1$ in $L_6$. By  Lemma~\ref{relations}(2),  $[(\r_0\r_1)^4,\r_2]=1$ and $[\r_0,(\r_0\r_1)^4]=1$ in $L_6$. Since $n-1-2-(n-4)=1$,   we have $|L_6|=2^{n-1}$ by taking $(n,s,t)=(n-1,2,n-4)$ for $H$ in Propositions~\ref{$2^s,2^t$}. It follows $|H_6|=|L_6|\cdot|\lg (\r_0\r_1) ^4\rg|\leq 2^{n-1}\cdot 2=2^n$.

Now we claim $|H_6|=2^n$ and the listed exponents are the true orders of the corresponding elements. To do that,  we will construct a permutation group $A$ of order at least $2^n$ that is an epimorphic image of $H_6$.

Set $t=2^{n-4}$ and write
$ci=\frac{t}{8}-i-1$ for $0\leq i\leq \frac{t}{8}-1$,
$i_{jt}^k=jt+8i+k$ and ${ci}_{jt}^k=jt+8ci+k$ for $0 \leq j \leq 3$ and $1 \leq k \leq 8$.  Note that $0\leq i\leq \frac{t}{8}-1$ if and only if $0\leq ci\leq \frac{t}{8}-1$. Clearly, $1\leq i_{jt}^k,{ci}_{jt}^k\leq 4t$. Let $A=\langle a,b,c\rangle$, where $a, b, c$ are permutations on the set $\{1,2,\cdots, 2^{n-2}\}$, defined as
\vskip 0.2cm
\begin{small}
$\begin{array}{rl}
  a=&\prod_{i=0}^{\frac{t}{8}-1}
    (i_{2t}^{1},ci_{3t}^{8})(i_{2t}^{8},ci_{3t}^{1})(i_{0}^{2},ci_{2t}^{7})
    (i_{t}^{2},i_{2t}^{2})(i_{3t}^{2},ci_{t}^{7})(i_{0}^{7},i_{3t}^{7})
    (i_{0}^{3},ci_{2t}^{6})(i_{t}^{3},i_{2t}^{3})\\
    &(i_{3t}^{3},ci_{t}^{6})(i_{0}^{6},i_{3t}^{6})(i_{0}^{4},ci_{t}^{5})(i_{0}^{5},ci_{t}^{4}),  \\ b=&\prod_{j=0}^3\prod_{i=0}^{\frac{t}{8}-1}(i_{jt}^1,i_{jt}^2)(i_{jt}^3,i_{jt}^4)
  (i_{jt}^5,i_{jt}^6)(i_{jt}^7,i_{jt}^8),          \\

c=&\prod_{j=0}^3[\prod_{i=0}^{\frac{t}{8}-1}(i_{jt}^2,i_{jt}^3)(i_{jt}^4,i_{jt}^5)
  (i_{jt}^6,i_{jt}^7) \cdot\prod_{i=0}^{\frac{t}{8}-2}(i_{jt}^8,(i+1)_{jt}^1)]\\
& \prod_{i=0}^{1} [(0_{0+2ti}^{1}) ((\frac{t}{8}-1)_{t+2ti}^{8})((\frac{t}{8}-1)_{0+2ti}^{8},0_{t+2ti}^{1})].\\

\end{array}$
\end{small}\\
Here, $(i+1)_{jt}^1=jt+8(i+1)+1$ for $0\leq i\leq \frac{t}{8}-2$.
It is easy to see that $a$ is fixed under conjugacy of $c$, that is, $a^c=a$. It follows $(ac)^2=1$.

Let $\a=a,b$, or $c$. Then $\a$ is an involution. Recall that $ci=\frac{t}{8}-i-1$. Since  $0\leq i\leq \frac{t}{8}-1$ if and only if $0\leq ci\leq \frac{t}{8}-1$, it is easy to see that if $\a$ interchanges $i_{j_1t}^{k_1}$ and $i_{j_2t}^{k_2}$ then $\a$ also interchanges $ci_{j_1t}^{k_1}$ and $ci_{j_2t}^{k_2}$, and if $\a$ interchanges $i_{j_1t}^{k_1}$ and $ci_{j_2t}^{k_2}$ then $\a$ also interchanges $ci_{j_1t}^{k_1}$ and $i_{j_2t}^{k_2}$. These facts are very helpful for the following computations.

\vskip 0.2cm
\begin{small}
$\begin{array}{lcl}
  ab&=&\prod_{i=0}^{\frac{t}{8}-1}
     (i_{0}^{1},i_{0}^{2},ci_{2t}^{8},i_{3t}^{2},ci_{t}^{8},ci_{t}^{7},i_{3t}^{1},ci_{2t}^{7})
     (i_{0}^{3},ci_{2t}^{5},ci_{2t}^{6},i_{0}^{4},ci_{t}^{6},i_{3t}^{4},i_{3t}^{3},ci_{t}^{5})\\

   & & (i_{0}^{5},ci_{t}^{3},ci_{2t}^{4},ci_{2t}^{3},ci_{t}^{4},i_{0}^{6},i_{3t}^{5},i_{3t}^{6})
     (i_{0}^{7},i_{3t}^{8},ci_{2t}^{2},ci_{t}^{1},ci_{t}^{2},ci_{2t}^{1},i_{3t}^{7},i_{0}^{8}),\\

bc&=&\prod_{j=0}^{1}(1+2tj,3+2tj,\cdots, 2t-1+2tj, 2t+2tj, 2t-2+2tj, \cdots, 2+2tj), \\

  (ab)^2&=&\prod_{i=0}^{\frac{t}{8}-1}
  (i_{0}^{1},ci_{2t}^{8},ci_{t}^{8},i_{3t}^{1})
  (i_{0}^{8},i_{3t}^{8},ci_{t}^{1},ci_{2t}^{1})
  (i_{0}^{2},i_{3t}^{2},ci_{t}^{7},ci_{2t}^{7})
  (i_{0}^{7},ci_{2t}^{2},ci_{t}^{2},i_{3t}^{7})\\

  &&(i_{0}^{3},ci_{2t}^{6},ci_{t}^{6},i_{3t}^{3})
  (i_{0}^{6},i_{3t}^{6},ci_{t}^{3},ci_{2t}^{3})
  (i_{0}^{4},i_{3t}^{4},ci_{t}^{5},ci_{2t}^{5})
  (i_{0}^{5},ci_{2t}^{4},ci_{t}^{4},i_{3t}^{5}),\\

  (ab)^4&=&\prod_{i=0}^{\frac{t}{8}-1}
  (i_{0}^{1},ci_{t}^{8},)
  (i_{t}^{1},ci_{0}^{8})
   (i_{2t}^{1},ci_{3t}^{8})
  (i_{3t}^{1},ci_{2t}^{8})
    (i_{0}^{2},ci_{t}^{7})
  (i_{t}^{2},ci_{0}^{7})
  (i_{2t}^{2},ci_{3t}^{7})
    (i_{3t}^{2},ci_{2t}^{7})
  \\

  &&(i_{0}^{3},ci_{t}^{6})
  (i_{t}^{3},ci_{0}^{6})
   (i_{2t}^{6},ci_{3t}^{3})
  (i_{3t}^{6},ci_{2t}^{3})
  (i_{0}^{4},ci_{t}^{5})
  (i_{t}^{4},ci_{0}^{5})
    (i_{2t}^{4},ci_{3t}^{5})
  (i_{3t}^{4},ci_{2t}^{5})
  ,\\

    ((ab)^2)^c&=&\prod_{i=0}^{\frac{t}{8}-1}
  (i_{0}^{1},i_{3t}^{1},ci_{t}^{8},ci_{2t}^{8})
  (i_{0}^{8},ci_{2t}^{1},ci_{t}^{1},i_{3t}^{8})
  (i_{0}^{2},ci_{2t}^{7},ci_{t}^{7},i_{3t}^{2})
  (i_{0}^{7},i_{3t}^{7},ci_{t}^{2},ci_{2t}^{2})
\\

  &&  (i_{0}^{3},i_{3t}^{3},ci_{t}^{6},ci_{2t}^{6})
  (i_{0}^{6},ci_{2t}^{3},ci_{t}^{3},i_{3t}^{6})
  (i_{0}^{4},ci_{2t}^{5},ci_{t}^{5},i_{3t}^{4})
  (i_{0}^{5},i_{3t}^{5},ci_{t}^{4},ci_{2t}^{4}).\\

\end{array}$
\end{small}
\vskip 0.2cm

The above computations imply $(ab)^8=1$ and $(bc)^{2^{n-3}}=1$. Furthermore, $(ab)^{-2}=c(ab)^2c$, that is, $[(ab)^2,c]=(ab)^{-4}=(ab)^4$. It is clear that $(bc)^{2^{n-4}}$ interchanges $i_0^k$ and $ci_{t}^{9-k}$ as $i_0^k+ci_{t}^{9-k}=2t+1$ (note that $1 \leq i_0^k\leq t$ and $t+1 \leq ci_{t}^{9-k}\leq 2t$), and similarly $(bc)^{2^{n-4}}$ interchanges $i_{2t}^{k}$ and $ci_{3t}^{9-k}$.
It is easy to check  $(bc)^{2^{n-4}}=(ab)^4=[(ab)^2,c]$.
So the generators $a,b,c$ of $A$ satisfy the same relations  as $\rho_0,\rho_1,\rho_2$ do in $H_6$, and hence there is an epimorphism $\a$ from $H_6$ to $A$ such that $\r_0^\a=a$, $\r_1^\a=b$ and $\r_2^\a=c$. Clearly, $A$ is transitive on $\{1,2,\cdots,2^{n-2}\}$ and $a,c\in A_1$, the stabilizer of $1$ in $A$. It follows that $|A|\geq 2^n$.  Since $|H_6|\leq 2^n$, $\a$ is an isomorphism and $|H_6|=|A|=2^n$.  Therefore the listed exponents are the true orders of the corresponding elements, as claimed.

Recall that $G_6=\lg \r_0, \r_1, \r_2 \ |\ R_2, (\r_0\r_1)^4(\r_1\r_2)^4(\r_1\r_2)^{2^{n-4}}, [(\r_0\r_1)^2,\r_2](\r_0\r_1)^4\rg$ with
$R_2=\{ \r_0^2, \r_1^2, \r_2^2, (\r_0\r_2)^2,(\r_0\r_1)^{2^{n-3}}, (\r_1\r_2)^{2^{n-3}}\}$. Note that in $H_6$, $(\r_1\r_0)^2((\r_0\r_1)^2)^{\r_2}=[(\r_0\r_1)^2, \r_2]=(\r_1\r_0)^4=(\r_0\r_1)^4=(\r_1\r_2)^{2^{n-4}}$. This implies $((\r_0\r_1)^2)^{\r_2}=(\r_0\r_1)^{-2}$ and $((\r_0\r_1)^2)^{\r_1\r_2}=(\r_0\r_1)^2$, that is, $[(\r_0\r_1)^2, \r_1\r_2]=1$. In particular, $[(\r_0\r_1)^2, (\r_1\r_2)^2]=1$. By Lemma~\ref{relations}(1), $[\r_0, (\r_1\r_2)^2]=[(\r_0\r_1)^2, \r_2]^{\r_1\r_2}=((\r_1\r_2)^{2^{n-4}})^{\r_1\r_2}=(\r_1\r_2)^{2^{n-4}}$, that is, $((\r_2\r_1)^2)^{\r_0}=(\r_1\r_2)^{2^{n-4}-2}$.

Now we have  $(\r_0\r_2\r_1)^2=((\r_2\r_1)^2)^{\r_0}(\r_0\r_1)^2=(\r_1\r_2)^{2^{n-4}-2}(\r_0\r_1)^2$ in $H_6$, and since $[(\r_0\r_1)^2, (\r_1\r_2)^2]=1$, it follows that $(\r_0\r_2\r_1)^4=(\r_1\r_2)^{-4}(\r_0\r_1)^4$ and $(\r_0\r_2\r_1)^8=(\r_1\r_2)^{-8}$
$(\r_0\r_1)^8=(\r_2\r_1)^8$ , which implies $o(\r_0\r_2\r_1)=o(\r_1\r_2)=2^{n-3}$ in $H_6$.
Thus,  $(\r_0\r_2\r_1)^4(\r_1\r_2)^4$
$(\r_1\r_2)^{2^{n-4}}=(\r_1\r_2)^{-4}(\r_0\r_1)^4(\r_1\r_2)^4(\r_1\r_2)^{2^{n-4}}=1$.
Since $\lg \r_0\r_1\rg\cap \lg\r_1\r_2\rg\unlhd H_6$, we have $\lg(\r_1\r_2)^{2^{n-4}}\rg\unlhd H_6$, and since  $\lg(\r_1\r_2)^{2^{n-4}}\rg\cong \mz_2$, we have $(\r_1\r_2)^{2^{n-4}}\in Z(H_6)$, the center of $H_6$. It follows
$[(\r_0\r_2\r_1)^2, \r_2](\r_0\r_2\r_1)^4=[(\r_2\r_1)^2(\r_0\r_1)^2, \r_2](\r_0\r_2\r_1)^4=[(\r_2\r_1)^2, \r_2]^{(\r_0\r_1)^2}$
$[(\r_0\r_1)^2, \r_2](\r_0\r_2\r_1)^4=(\r_1\r_2)^4(\r_0\r_1)^4(\r_1\r_2)^{-4}(\r_0\r_1)^4=1$.
So the generators $\r_0\r_2, \r_1, \r_2$ in $H_6$ satisfy the same relations as  $\r_0, \r_1, \r_2$ do in $G_6$, and there is an epimorphism $\phi: G_6 \mapsto H_6$ such that $\r_0^{\phi}=\r_0\r_2$, $\r_1^{\phi}=\r_1$ and $\r_2^{\phi}=\r_2$. Since $|H_6|=2^n$ and $|G_6| \leq 2^n$,  $\phi$ is an isomorphism. It follows $|G_6|=2^n$, and since $(\r_0\r_1)^\phi=\r_0\r_2\r_1$, we have $o(\r_0\r_1)=2^{n-3}$ and $o(\r_0\r_1)=2^{n-3}$; furthermore the listed exponents in $G_6$ are the true orders of the corresponding elements. This finishes the proof of  sufficiency part in (2).

To prove the necessity  part in (2), let $G$ be the automorphism group of a regular map of order $2^{n}$ and  type $\{2^{n-3}, 2^{n-3}\}$. Then $o(\r_0)=o(\r_1)=o(\r_2)=o(\r_0\r_2)=2$ and $o(\r_0\r_1)=o(\r_1\r_2)=2^{n-3}$. We only need to show $G=G_3, G_4, G_5$ or $G_6$, and it will be done by induction on $n$. This is true for $n=10$ by {\sc Magma}.

Assume $n\geq 11$. Take $N=\lg (\r_1\r_2)^{2^{n-4}}\rg$. By Lemma~\ref{quotient}, $N\unlhd G$ and $\olg=G/N$ is the automorphism group of a regular map of order $2^{n-1}$ with type $\{2^{n-4}, 2^{n-4}\}$. Since $|\olg|=2^{n-1}$, by induction hypothesis we may assume $\olg=\olg_3, \olg_4, \olg_5$ or $\olg_6$ with $R=\{\overline{\r_0}^2, \overline{\r_1}^2, \overline{\r_2}^2, (\overline{\r_0}\overline{\r_1})^{2^{n-4}}, (\overline{\r_1}\overline{\r_2})^{2^{n-4}}, (\overline{\r_0}\overline{\r_2})^2\}$, where

\begin{itemize}
  \item [$\olg_3$]$=\lg \overline{\r_0}, \overline{\r_1}, \overline{\r_2} \ |\ R, (\overline{\r_0}\overline{\r_1})^4(\overline{\r_1}\overline{\r_2})^4, [(\overline{\r_0}\overline{\r_1})^2, \overline{\r_2}](\overline{\r_0}\overline{\r_1})^4\rg$,
 \item [$\olg_4$]$=\lg \overline{\r_0}, \overline{\r_1}, \overline{\r_2} \ |\ R, (\overline{\r_0}\overline{\r_1})^4(\overline{\r_1}\overline{\r_2})^4,  [(\overline{\r_0}\overline{\r_1})^2, \overline{\r_2}](\overline{\r_0}\overline{\r_1})^4(\overline{\r_1}\overline{\r_2})^{2^{n-5}}\rg$,
   \item [$\olg_5$]$=\lg \overline{\r_0}, \overline{\r_1}, \overline{\r_2} \ |\ R, (\overline{\r_0}\overline{\r_1})^4(\overline{\r_1}\overline{\r_2})^4(\overline{\r_1}\overline{\r_2})^{2^{n-5}}, [(\overline{\r_0}\overline{\r_1})^2, \overline{\r_2}](\overline{\r_0}\overline{\r_1})^4(\overline{\r_1}\overline{\r_2})^{2^{n-5}}\rg$,

   \item [$\olg_6$]$=\lg \overline{\r_0}, \overline{\r_1}, \overline{\r_2} \ |\ R, (\overline{\r_0}\overline{\r_1})^4(\overline{\r_1}\overline{\r_2})^4(\overline{\r_1}\overline{\r_2})^{2^{n-5}},  [(\overline{\r_0}\overline{\r_1})^2, \overline{\r_2}](\overline{\r_0}\overline{\r_1})^4\rg$.
\end{itemize}

Let $H=\lg \r_0\r_1, \r_1\r_2\rg$ be the rotation subgroup of $G$. Then $|G:H| \leq 2$. Note that $|\lg \r_0, \r_1\rg|=|\lg \r_1, \r_2\rg|=2^{n-2}$, $|\lg \r_0, \r_2\rg|=4$ and $|G|=2^n$. This implies that $\{\r_0,\r_1,\r_2\}$ is a minimal generating set of $G$. By  Proposition~\ref{burnside}~(2), $G$ has rank $3$, and since $H$ is generated by two elements, we have $|G:H|=2$, that is, $|H| = 2^{n-1}$. Since $|\lg \r_0\r_1\rg \lg \r_1\r_2\rg|=\frac{|\lg \r_0\r_1\rg||\lg \r_1\r_2\rg|}{|\lg \r_0\r_1\rg \cap \lg \r_1\r_2\rg|}=\frac{2^{2n-6}}{|\lg \r_0\r_1\rg \cap \lg \r_1\r_2\rg|} \leq |H|=2^{n-1}$, we have $|\lg \r_0\r_1\rg \cap \lg \r_1\r_2\rg| \geq 2^{n-5}$. Then $\lg \r_0\r_1\rg \cap \lg \r_1\r_2\rg$ has a subgroup of order $2^{n-5}$, which is the unique subgroup of order $2^{n-5}$ in $\lg \r_0\r_1\rg$ and $\lg \r_1\r_2\rg$ respectively, that is, $\lg (\r_0\r_1)^4\rg$ and $\lg (\r_1\r_2)^4\rg$. It follows $\lg (\r_0\r_1) ^4\rg =\lg (\r_1\r_2)^4\rg$.

Suppose $\olg=\olg_4$ or $\olg_5$. Then $[(\r_0\r_1)^2, \r_2](\r_0\r_1)^4(\r_1\r_2)^{2^{n-5}}=1$ or $(\r_1\r_2)^{2^{n-4}}$
because $\lg (\r_1\r_2)^{2^{n-4}}\rg \cong \mathbb{Z}_2$. It follows $[(\r_0\r_1)^2, \r_2](\r_0\r_1)^4=(\r_1\r_2)^{\d \cdot 2^{n-5}}$ with $\d=1$ or $-1$. Since $\lg (\r_0\r_1) ^4\rg =\lg (\r_1\r_2)^4\rg$ and $[(\r_0\r_1)^2, \r_2]\in \lg (\r_0\r_1) ^4\rg$, we have $[(\r_0\r_1)^2, \r_2]^{\r_0}=[(\r_0\r_1)^2, \r_2]^{\r_1}=[(\r_0\r_1)^2, \r_2]^{\r_2}=[(\r_0\r_1)^2, \r_2]^{-1}$.
By Proposition~\ref{commutator}, $[(\r_0\r_1)^2, (\r_1\r_2)^2]=[(\r_0\r_1)^2, \r_2][(\r_0\r_1)^2, \r_1\r_2\r_1]^{\r_2}=[(\r_0\r_1)^2, \r_2][(\r_0\r_1)^2, \r_2]^{\r_0\r_1\r_2}=[(\r_0\r_1)^2, \r_2][(\r_0\r_1)^2, \r_2]^{-1}=1$. On the other hand, $[(\r_0\r_1)^2, \r_1]=(\r_1\r_0)^4\in \lg (\r_1\r_2)^4\rg$, and hence $[(\r_0\r_1)^2, \r_1]^{\r_2}=(\r_1\r_0)^{-4}$ and
$[(\r_0\r_1)^2, \r_1\r_2]=[(\r_0\r_1)^2, \r_2][(\r_0\r_1)^2, \r_1]^{\r_2}=(\r_1\r_2)^{\d \cdot 2^{n-5}}(\r_1\r_0)^4(\r_1\r_0)^{-4}=(\r_1\r_2)^{\d \cdot 2^{n-5}}$. It follows $1=[(\r_0\r_1)^2, (\r_1\r_2)^2]=[(\r_0\r_1)^2, \r_1\r_2][(\r_0\r_1)^2, \r_1\r_2]^{\r_1\r_2}=(\r_1\r_2)^{2^{n-4}}$, which is impossible because $o(\r_1\r_2)=2^{n-3}$.

Suppose $\olg=\olg_6$. Since $\lg (\r_1\r_2)^{2^{n-4}}\rg \cong \mathbb{Z}_2$, we have $(\r_0\r_1)^4(\r_1\r_2)^4(\r_1\r_2)^{2^{n-5}}=1$ or $(\r_1\r_2)^{2^{n-4}}$, and $[(\r_0\r_1)^2, \r_2](\r_0\r_1)^4=1$ or $(\r_1\r_2)^{2^{n-4}}$,
which implies $(\r_0\r_1)^4(\r_1\r_2)^4=(\r_1\r_2)^{\d \cdot 2^{n-5}}$ with $\d=1$ or $-1$ and $[(\r_0\r_1)^2, \r_2]=(\r_1\r_2)^{\gamma\cdot 2^{n-4}}(\r_1\r_0)^4$ with $\gamma=0$ or $1$. Since $\lg(\r_0\r_1)^4\rg=\lg(\r_1\r_2)^4\rg$, we have $[(\r_0\r_1)^2, \r_2]\in \lg(\r_1\r_2)^4\rg$. By  Lemma~\ref{relations}(1), we have $[\r_0, (\r_1\r_2)^2]=[(\r_0\r_1)^2, \r_2]^{\r_1\r_2}=[(\r_0\r_1)^2, \r_2]$.

It follows, by
Proposition~\ref{commutator}, that $[(\r_0\r_1)^2, (\r_1\r_2)^2]=[(\r_0\r_1)^2, \r_2][(\r_0\r_1)^2, \r_2]^{\r_0\r_1\r_2}=[(\r_0\r_1)^2, \r_2][(\r_0\r_1)^2, \r_2]^{-1}=1$, and also $[\r_0, (\r_1\r_2)^4]
=[\r_0, (\r_1\r_2)^2][\r_0, (\r_1\r_2)^2]^{(\r_1\r_2)^2}
=[(\r_0\r_1)^2, \r_2]^2=(\r_1\r_0)^8$. Since $[\r_0, (\r_1\r_2)^4]=((\r_2\r_1)^4)^{\r_0}(\r_1\r_2)^4
=(\r_2\r_1)^{-4}(\r_1\r_2)^4=(\r_1\r_2)^8$ as $\lg(\r_0\r_1)^4\rg=\lg(\r_1\r_2)^4\rg$, we have  $1=(\r_0\r_1)^8(\r_1\r_2)^8=
((\r_0\r_1)^4(\r_1\r_2)^4)^2=((\r_2\r_1)^{\d \cdot 2^{n-5}})^2=(\r_2\r_1)^{2^{n-4}}$, which is
impossible because $o(\r_1\r_2)=2^{n-3}$.

Thus, $\olg=\olg_3$. Since $\langle (\r_1\r_2)^{2^{n-4}}\rangle\cong\mz_2$, we have $[(\r_0\r_1)^2, \r_2](\r_0\r_1)^4=1$ or $(\r_1\r_2)^{2^{n-4}}$ and $(\r_0\r_1)^4(\r_1\r_2)^4=1$ or $(\r_1\r_2)^{2^{n-4}}$, and hence $G=G_3, G_4, G_5$ or $G_6$.  \hfill\qed

\medskip
\f {\bf Acknowledgements:} This work was partially supported by the National Natural Science Foundation of China (11571035, 11731002), the 111 Project of China (B16002), and the third author was supported by Basic Science Research
Program through the National Research Foundation of Korea
funded by the Ministry of Education (2015R1D1A1A09059016).


\begin{thebibliography}{99}



\bibitem{DCJ2014}
Archdeacon, D., Bonnington, C.P., \v{S}ir\'{a}\v{n}, J.:
Regular pinched maps.
Australas. J. Combin. {\bf 58}, 16--26 (2014)

\bibitem{BDR}
Ban, Y.F., Du, S.F., Liu, Y., Nedela, R., \"{S}koviera, M.: Classification of regular maps whose automorphism groups are $2$-groups of class three, in preparation.

\bibitem{GroupBookss} Berkovich, Y.:
 Groups of Prime Power Order,
 vol. 1.Walter de Gruyter, Berlin (2008)


\bibitem{BCP97}
Bosma, W., Cannon, J., Playoust, C.:
 The {M}agma {A}lgebra {S}ystem. {I}:  the user language.
 J. Symb. Comput. {\bf 24}, 235--265 (1997)


\bibitem{H1927}
Brahana, H.:
Regular maps and their groups.
Amer. J. Math. {\bf 49}, 268--284 (1927).

\bibitem{RD1985}
Bryant, R.P., Singerman, D.: Foundations of the theory of maps on surfaces with boundary.
Q. J. Math. {\bf 36}, 17--41 (1985)

\bibitem{ARJ2005}
Breda d'Azevedo, A., Nedela, R., \v{S}ir\'{a}\v{n}, J.: Classification of regular maps of negative prime Euler characteristic.
Trans. Amer. Math. Soc. {\bf 357}, 4175--4190 (2005)


\bibitem{atles}Conder, M.D.E.: Regular maps (on orientable or non-orientable surfaces) with up to 1000 edges, available at \url{https://www.math.auckland.ac.nz/~conder/RegularMapsWithUpTo1000Edges.txt}.



\bibitem{CMa2015}
Conder, M.D.E., Ma, J.C.: Regular maps with simple underlying graphs.
J. Combin. Theory Ser. B {\bf 110}, 1--18 (2015)


\bibitem{CDNS}
Conder, M.D.E., Du, S.F., Nedela, R., \v{S}koviera, M.: Regular maps with nilpotent automorphism group.
J. Algebraic. Combin. {\bf 44}, 863--874 (2016)


\bibitem{MRJ2012}
Conder, M.D.E.,  Nedela, R., \v{S}ir\'{a}\v{n}, J.: Classification of regular maps of Euler characteristic -3p.
J. Combin. Theory Ser. B {\bf 102}, 967--981 (2012)

\bibitem{CHNS2012}
Conder, M.D.E., Huc\'{i}kov\'{a}, V., Nedela, R., \v{S}ir\'{a}\v{n}, J.: Chiral maps of given hyperbolic type, Bull. Lond. Math. Soc. {\bf 48}, 38--52 (2016)


\bibitem{MP2001}
Conder, M.D.E., Dobcs\'{a}nyi, P.: Determination of all regular maps of small genus.
J. Combin. Theory Ser. B {\bf 81}, 224--242 (2001)



\bibitem{2GroupBookss} Doerk, K., Hawkes, T.:
Finite Soluble Groups,
Walter de Gruyter, Berlin (1992)






\bibitem{ARJM1999}
Gardiner, A., Nedela, R., \v{S}ir\'{a}\v{n}, J., \v{S}koviera, M.: Characterization of graphs which underlie regular maps on closed surfaces.
J. Lond. Math. Soc. {\bf 59}, 100--108 (1999)


\bibitem{NG2013}
Gill, N.: Orientably regular maps with Euler characteristic divisible by few primes.
J. Lond. Math. Soc. {\bf 88}, 118--136 (2013)

\bibitem{GW}
Gray, A., Wilson, S.: A more elementary proof of Gr\"{u}nbaum¡¯s conjecture. Congr. Numer. {\bf 72}, 25--32 (1990)

\bibitem{HFL}
Hou, D.-D., Feng, Y.-Q., Leemans, D.: Existence of regular 3-polytopes of order $2^n$. J. Group Theory, accepted. (2018)

\bibitem{HW}
Hu, K., Wang, N.: Classification of regular maps whose automorphism groups are $2$-groups of maximal class.
Acta Univ. M. Belii Ser. Math. {\bf 20}, 11--17 (2012)

\bibitem{GroupBooks} Huppert, B.:
 Endliche Gruppen \uppercase\expandafter{\romannumeral1},
 Springer, Berlin (1967)

\bibitem{G1994}
Jones, G.A.: Ree groups and Riemann surfaces.
J. Algebra {\bf 165}, 41--62 (1994)

\bibitem{JS1} Jones, G.A., Singerman, D.:
             Theory of maps on orientable surfaces.
           Proc. Lond. Math. Soc.  {\bf 37}, 273--307 (1978)


\bibitem{JNS}
Jendrol¡¯, S., Nedela, R., \v{S}koviera, M.: Constructing regular maps and graphs from planar quotients. Math. Slovaca. {\bf 47}, 155--170 (1997)







\bibitem{LCH}
Li, C.H., \v{S}ir\'{a}\v{n}, J.: Regular maps whose groups do not act faithfully on vertices, edges, or faces.
European J. Combin. {\bf 26}, 521--541 (2005)


\bibitem{ARM2012}
Malni\v{c}, A., Nedela, R., \v{S}koviera, M.: Regular maps with nilpotent automorphism groups.
European J. Combin. {\bf 33}, 1974--1986 (2012)


\bibitem{N} Nedela, R.: Regular maps - combinatorial objects relating  different fields of mathematics. J. Korean Math. Soc.  {\bf 38}, 1069--1105 (2001)



\bibitem{C1969}
Sah, Ch.-H.: Groups related to compact Riemann surfaces.
Acta Math. {\bf 123}, 13--42 (1969)

\bibitem{Si}
 \v{S}ir\'{a}\v{n}, J.: Non-orientable Regular Maps of a Given Type over Linear Fractional Groups.
Graphs Combin. {\bf 26}, 597--602 (2010)

\bibitem{Vi1983}
Vince, A.: Regular combinatorial maps. J. Combin. Theory Ser. B {\bf 35}, 256--277 (1983)

\end{thebibliography}
\end{document}